\documentclass{amsart}
\usepackage{pb-diagram}
\usepackage{amssymb}

\def\A{\mathcal {A}}
\def\B{\mathcal {B}}

\def\Fe{\mathcal {F}}
\def\R{\mathbf{R}}
\def\C{\mathbf{C}}
\def\F{\mathbf{F}}

\def\N{\mathbf{N}}
\def\Z{\mathbf{Z}}
\def\w{\varpi}
\def\p{\partial}
\def\pp#1#2{\frac{\p #1}{\p #2}}
\def\PB{\left\{\cdot\,,\cdot\right\}}
\def\pb#1{\left\{#1\right\}}
\def\LB{\left[\cdot\,,\cdot\right]}
\def\lb#1{\left[#1\right]}
\def\ideal#1{\langle #1\rangle}
\def\Cas{\mathop{\rm Cas}\nolimits}
\def\Char{\mathop{\rm char}\nolimits}
\def\Hom{\mathop{\rm Hom}\nolimits}
\def\div{\mathop{\rm Div}\nolimits}
\def\Vect{{\mathfrak{X}}}

\def\diff{{{\rm d}}}
\def\de{\delta}
\def\L{{\mathcal L}}
\def\v{\varphi}
\def\vn{\vec{\nabla}}

\newtheorem{thm}{Theorem}[section]
\newtheorem{prp}[thm]{Proposition}
\newtheorem{cor}[thm]{Corollary}
\newtheorem{lma}[thm]{Lemma}
\theoremstyle{definition}
\newtheorem{dfn}[thm]{Definition}
\newtheorem{rem}[thm]{Remark}

\newenvironment{eqn*}[1][1.5]
  {$$\renewcommand{\arraystretch}{#1}
      \begin{array}{rcl}}
      {\end{array}$$}

\newenvironment{eqn}[2][1.5]
  {\begin{equation}\label{#2}
   \renewcommand{\arraystretch}{#1}
   \begin{array}{rcl}}
  {\end{array}\end{equation}}

\begin{document}
\nocite{*}

\title{Poisson (co)homology and isolated singularities}
\author{Anne Pichereau}
  \address{Universit\'e de Poitiers,
           Laboratoire de Math\'ematiques et Applications,
           UMR 6086,
           T\'el\'eport 2,
           Boulevard Marie et Pierre Curie,
           BP 30 179,
           F-86 962 Futuroscope Chasseneuil Cedex}
  \email{anne.pichereau@math.univ-poitiers.fr}
\keywords{Poisson cohomology, Poisson homology, isolated singularities}
\subjclass[2000]{17B55, 17B63}
\begin{abstract}
To each polynomial $\v\in\F[x,y,z]$ is associated a Poisson structure on
$\F^3$, a surface and a Poisson structure on this surface. When
$\v$ is weight homogeneous with an isolated singularity, we determine the
Poisson cohomology and homology of the two Poisson varieties obtained.
\end{abstract}
\maketitle
\tableofcontents

\section{Introduction}
\label{intro}

The first Poisson structures appeared in classical mechanics.
In $1809$, D. Poisson introduced a bracket
of functions, given by:
\begin{eqnarray}\label{first_poisson}
\pb{f,g}=\sum_{i=1}^r \left(\pp{f}{q_i}\pp{g}{p_i}-\pp{f}{p_i}\pp{g}{q_i}\right),
\end{eqnarray}
for two smooth functions $f,g$ on $\R^{2r}$. It permits one to write the
Hamilton's equations as differential equations, where positions
($q_i$) and impulsions ($p_i$) play symmetric roles. Indeed,
denoting by $H$ the total energy of the system, these equations become:
\begin{eqnarray*}
 \begin{array}{cc}
 \renewcommand{\arraystretch}{1.2}
\begin{array}{ccl}   
  \dot{q_i} &=& \pb{q_i,H},\\
  \dot{p_i} &=& \pb{p_i,H}, 
 \end{array}&
\quad 1\leq i\leq r.
 \end{array}
\end{eqnarray*}
D. Poisson also pointed out that if $f$ and $g$ are constants of motion,
then $\pb{f,g}$ is also a constant of motion and this
phenomenon was explained in $1839$ by C. Jacobi, who proved that 
(\ref{first_poisson}) satisfies what is now called the Jacobi identity:
\begin{eqnarray}\label{intro_jacobiator}
\pb{\pb{f,g},h}+\pb{\pb{g,h},f}+\pb{\pb{h,f},g}=0.
\end{eqnarray}
This important identity leads to the definition of a Poisson algebra
as an algebra $\B$ equipped with a skew-symmetric biderivation $\PB$,
satisfying (\ref{intro_jacobiator}), for
all $f,g,h$, elements of $\B$. Said differently, a Poisson algebra is a Lie
algebra $(\B,\PB)$, where $\PB$ satifies the Leibniz rule 
$\pb{fg,h}=f\pb{g,h}+\pb{f,h}g$, for all $f,g,h\in\B$.
One talks about a Poisson variety, when its
algebra of functions is equipped with a Poisson structure. This notion
generalizes the notion of symplectic manifold.

For a given Poisson algebra $(\B,\PB)$, one defines a cohomology,
called Poisson cohomology, introduced by A. Lichnerowicz in
\cite{lichne}; see also \cite{huebs} for an algebraic approach. 
The cochains are the skew-symmetric multiderivations of
$\A$ and the coboundary operator is $-\lb{\pi,\,\cdot}_S$, where
$\pi:=\PB$ is the Poisson bracket and $\LB_S$ is the Schouten
bracket.
The resulting Poisson complex, defined in detail in Section \ref{ourstruct}, can be
viewed as the contravariant version of the de Rham complex. Its
cohomology gives very interesting information about the Poisson
structure, as for small $k$, the $k$-th Poisson cohomology space $H^k(\B,\pi)$
has the following interpretation:
\begin{eqnarray*}
H^0(\B,\pi) &=& \{\hbox{Casimir functions}\} := \{f\in\B\mid\pb{f,\,\cdot\ }=0\},\\
{}\\
H^1(\B,\pi) &=& \frac{\{\hbox{Poisson derivations}\}}{\{\hbox{Hamiltonian
                 derivations}\}},\\
{}\\
H^2(\B,\pi) &=& \frac{\{\hbox{skew-symmetric biderivations compatible with }\pi\}}
                  {\{\hbox{Lie derivatives of }\pi\}},\\
{}\\
H^3(\B,\pi) &=& \{\hbox{Obstructions to deformations of Poisson
                 structures}\}.
\end{eqnarray*}
\smallskip
Moreover, $H^2(\B,\pi)$ is fundamental in the study of normal forms of
Poisson structures (see \cite{zung}). We also denote by $\Cas(\B,\pi)$
the space of all Casimir functions of $(\B,\PB)$  (that is to say
$H^0(\B,\pi)$) and we point out that each $H^k(\B,\pi)$ is a
$\Cas(\B,\pi)$-module in a natural way.

To determine the Poisson cohomology of a given
Poisson algebra explicitly is, in general, difficult. One of the reasons
seems to be that Poisson cohomology is not a functor: a morphism 
$\pi:\A_1\rightarrow\A_2$ between Poisson algebras does not lead to a
morphism between their cochains (multiderivations), nor between their
corresponding Poisson cohomology groups.
In a few specific cases, Poisson cohomology has been determined. 
For a symplectic manifold, there exists a natural isomorphism between Poisson and de Rham
cohomology (see \cite{lichne}). In \cite{vais} and \cite{xu2}, one finds
some partial results about the case of regular
Poisson manifolds, while, for
Poisson-Lie groups, one can refer to \cite{gw}. Finally, the Poisson
cohomology in dimension two was computed in the germified and algebraic
cases in \cite{monnier} and \cite{r_v}.

Our purpose is to determine the Poisson cohomology of two classes of Poisson
varieties, intimately linked. The first class is composed of the singular surfaces
$\Fe_{\v}:\{\v=0\}$ in $\F^3$ ($\F$ is a field of
characteristic zero) that are 
defined by the zeros of polynomials $\v\in\F[x,y,z]$ and the second one
is the class of the Poisson varieties that are the ambient space $\F^3$, 
equipped with Poisson structures associated
to each $\v$. It means that we consider Poisson structures on the
algebras of regular functions on $\Fe_{\v}$ and $\F^3$, given by 
$\A_{\v}:=\F[x,y,z]/\ideal{\v}$ and  $\A:=\F[x,y,z]$ and that we
determine the Poisson cohomology of the Poisson algebras obtained.

We point out that the dimension three is the first one in which there
is a real condition for a biderivation to be a Poisson biderivation. The
Jacobi identity is indeed trivial in dimension two and every 
polynomial $\psi\in\F[x,y]$ leads to a Poisson
structure on the affine
space $\F[x,y]$, given by $\psi\pp{}{x}\wedge\pp{}{y}$. One
can consider the singular locus of such a structure, 
given by $\Gamma_{\psi}:\{\psi=0\}$. 
In \cite{r_v}, the authors determine the dimensions of the Poisson 
cohomology spaces, when $\psi$ is a homogeneous polynomial. They observe
that these dimensions are linked to the type of the
singularity of $\Gamma_{\psi}$. Conversely, in our context, we consider a
surface $\Fe_{\v}$, with a singularity, and a Poisson bracket
that do not bring other singularities. That is to say, this Poisson structure is
symplectic everywhere except on the singularities of $\Fe_{\v}$. In fact, it will
be the restriction of a Poisson structure $\PB_{\v}$ on $\F^3$, which is 
completely defined by the brackets:
\begin{eqnarray}
\label{brackets_chi_psi}
\pb{x,y}_{\v}=\pp{\varphi}{z},\quad \pb{y,z}_{\v}=\pp{\varphi}{x},
\quad \pb{z,x}_{\v}=\pp{\varphi}{y}, \quad (\v\in\A).
\end{eqnarray}
We suppose that $\Fe_{\v}$ has only one weight homogeneous
isolated singularity (at the origin). In fact, the hypothesis is that
$\v$ is a weight homogeneous polynomial with an isolated singularity. 

An other way to approach our context is to consider the Poisson structures
on $\A$ that admit a weight homogeneous Casimir and a singular locus
reduced to the origin. That leads to study the Poisson structures of
the form $\PB_{\v}$, with $\v$ weight homogeneous with an isolated
singularity.
As $\v$ is a Casimir for this structure, $\ideal{\v}$ is a Poisson
ideal of the Poisson algebra $(\A,\PB_{\v})$. This implies that
$\PB_{\v}$ goes down to the quotient algebra
$\A_{\v}=\F[x,y,z]/\ideal{\v}$. The singular surface $\Fe_{\v}$ is then 
the union of a symplectic leave of $\PB_{\v}$ and the origin.

For each $\v\in\A$ weight homogeneous with an isolated singularity, 
what we determine is the Poisson cohomology of both the Poisson
algebras introduced. Moreover, we turn these results to good account
to give the Poisson homology of these algebras. 
The Poisson cohomology spaces are respectively
denoted by $H^k(\A,\v)$ for $(\A,\PB_{\v})$ and $H^k(\A_{\v})$ for
the singular surface, while the Poisson homology
spaces are denoted by $H_k(\A,\v)$ and $H_k(\A_{\v})$. 

To develop a first idea about our results, one may think of $\v$ as a
homogeneous polynomial, of degree denoted by $\w(\v)$, such that 
its three partial derivatives have
only one common zero that is the
origin. This implies that
$$
\A_{sing}:=\A/\ideal{\pp{\v}{x},\pp{\v}{y},\pp{\v}{z}}
$$
is a finite dimensional $\F$-vector space. Its dimension is the
so-called Milnor number~$\mu$ (see \cite{miln}). This space gives information about the
(isolated) singularity of the surface $\Fe_{\v}$ (like multiplicity,
see also \cite{C_L_OS}) as it is exactly the algebra of regular functions on
this singularity. It plays an important role in the
Poisson cohomology of the algebra $(\A,\PB_{\v})$, so that this Poisson
cohomology is closely related to the type of the singularity of $\Fe_{\v}$.
We consider a family $u_0=1,u_1,\dots,u_{\mu-1}$ of homogeneous elements of $\A$,
whose images in $\A_{sing}$ give a $\F$-basis
of this $\F$-vector space.

The algebra of Casimir functions of the algebra $(\A,\PB_{\v})$ is given in
Proposition \ref{calculh0} and is simply the
algebra generated by $\v$, that is to say 
$\Cas(\A,\v)=H^0(\A,\v)\simeq\bigoplus_{i\in\N}\F\v^i$. 
In Proposition \ref{calculh1}, we see that the first Poisson
cohomology space of $\A$ is equal to zero if the
degree of $\v$, $\w(\v)$, is equal to $3$ and otherwise $H^1(\A,\v)$
is the $\Cas(\A,\v)$-module given by 
$$
H^1(\A,\v)\simeq\Cas(\A,\v)\vec{e},
$$
where $\vec{e}:=(x,y,z)$ corresponds to the Euler derivation 
$x\pp{}{x}+y\pp{}{y}+z\pp{}{z}$. 
Notice that the cubic polynomials play a special role here; in the weight
homogeneous case, this role is played by the polynomials of degree the sum
of the weights of the three variables $x,y,z$.
Moreover, with
Proposition \ref{calculh2}, we see that the case $\w(\v)=3$ is also
the unique case where the biderivation $\PB_{\v}$ is not an exact
Poisson structure, i.e.\ $\PB_{\v}$, which is a
$2$-cocycle of the Poisson cohomology of $(\A,\PB_{\v})$,
is not a $2$-coboundary (see \cite{huebs}). 
Proposition \ref{calculh2} affirms indeed that the second Poisson
cohomology space is exactly 
\begin{eqnarray*}
H^2(\A,\v) &\simeq& \bigoplus_{\stackrel{j\geq 1}{\w(u_j)\not= \w(\v)-3}} 
    \Cas(\A,\v)\vn u_j \oplus  \bigoplus_{\w(u_j)=\w(\v)-3} 
    \Cas(\A,\v)u_j\vn\v\\
&&\qquad\qquad\qquad\qquad 
\oplus  \bigoplus_{\stackrel{j\geq 1}{\w(u_j)=\w(\v)-3}}\F\vn u_j.
\end{eqnarray*}
This writing has been obtained from the third Poisson cohomology
space, which is determined in Proposition
\ref{calculh3}, and is exactly the free $\Cas(\A,\v)$-module
$$
H^3(\A,\v)\simeq \Cas(\A,\v)\otimes_{\F}\A_{sing}.
$$
It may be remarked that $H^2(\A,\v)$ is the unique Poisson cohomology
space of $(\A,\PB_{\v})$ which is not always a free module over the algebra of Casimirs. 

In Chapter \ref{surface}, we give the Poisson cohomology spaces of
the singular surface $\Fe_{\v}$, by considering the algebra $\A_{\v}$. 
For this Poisson algebra, the
Casimirs are simply the elements of $\F$ and, according to Propositions
\ref{calculh1surf} and \ref{calculh2surf}, we have:
$$
H^1(\A_{\v})\simeq\bigoplus_{\w(u_j)=\w(\v)-3}
            \F u_j\,\vec{e}\;, \qquad
H^2(\A_{\v})\simeq\bigoplus_{\w(u_j)=\w(\v)-3} \F u_j \vn\v.
$$
Finally, in Chapter \ref{homology}, we determine the Poisson homology of the
algebra $(\F^3,\PB_{\v})$ and of the singular surface $\Fe_{\v}$. 
We explain first, in Proposition \ref{homology_A}, that we have isomorphisms
$$
H_k(\A,\v)\simeq H^{3-k}(\A,\v), \hbox{ for all } k=0,1,2,3.
$$
Then, using the results about Poisson cohomology of $(\A,\PB_{\v})$, 
we compute the Poisson homology spaces of $\Fe_{\v}$ and we obtain, in
Proposition \ref{homology_surf},
$$
H_0(\A_{\v}) \simeq H_2(\A_{\v})\simeq \A_{sing}\;;\qquad
H_1(\A_{\v}) \simeq \bigoplus_{j=1}^{\mu-1}\F\,\vn u_j.
$$

Since the coboundary operator is
a weight homogeneous operator (see Section \ref{homder}), all
our arguments remain true if we replace the algebra $\A=\F[x,y,z]$ by
the algebra of all formal power series $\bar{\A}:=\F[[x,y,z]]$, still equipped with
the Poisson structure $\PB_{\v}$, with $\v$ a weight homogeneous element of
$\A$. 
It suffices to replace $\Cas(\A,\v)=\F[\v]$ by
$\Cas(\bar{\A},\v)=\F[[\v]]$, the algebra of formal power series in $\v$.\\

I would like to take the opportunity to thank my thesis advisor, 
Pol Vanhaecke, for suggesting to me this interesting problem and
for his availability all along this project. 
I am also indebted to Claude
Quitt\'e, whose knowledge of regular sequences was precious for me, and
Camille Laurent for his explanations about the modular class.

I finally would like to thank Prof. M. van den Bergh. After writing this
paper, he pointed out to me that, in his article ``Noncommutative homology
of some three-dimensional quantum spaces'' (see \cite{vdB}), he computed the Poisson
homology spaces of the Poisson algebra $(\A,\PB_{\v})$, for 
$\varphi=\frac{q_1}{3}(x^3+y^3+z^3)+2p_1xyz$, where $p_1$ and $q_1$ are
parameters. This case is a particular one
of the Poisson homology that I determine, and the method is very similar.
\section{The Poisson cohomology complex associated to a polynomial}
\label{Poissoncomplex}
\subsection{Poisson structures on $\A=\F[x,y,z]$ and their cohomology}
\label{ourstruct}

Let $\A$ be the polynomial algebra $\A=\F[x,y,z]$, where $\F$ is a
field of characteristic zero  and let $\v\in\A$. A Poisson structure
on $\A$ is defined by the brackets: 
\begin{eqnarray}
\label{bracket}
\pb{x,y}_{\v}=\pp{\varphi}{z},\quad \pb{y,z}_{\v}=\pp{\varphi}{x},
\quad \pb{z,x}_{\v}=\pp{\varphi}{y}.
\end{eqnarray}
Recall that a Poisson bracket on an associative and commutative algebra $\B$ is
a skew-symmetric bilinear map
$\PB$, from $\B^2$ to $\B$ (element of $\Hom(\wedge^2\B,\B)$), which is a
derivation in each of its arguments and which satisfies the Jacobi identity: 
\begin{eqnarray}
\label{jacobi}
\pb{\pb{f,g},h}+\pb{\pb{g,h},f}+\pb{\pb{h,f},g}=0,
\end{eqnarray}
for each $f,g,h\in\B$. In the particular case of $\A$, 
the brackets of the generators $x, y, z$ define totally the Poisson
bracket, in view of the derivation property, and moreover 
the Jacobi identity is satisfied for all 
$f,g,h\in\A$ if and only if it is satisfied for $x, y, z$ (see \cite{pol}). Here, 
an easy computation shows that this condition is satisfied by the
bracket $\PB_{\v}$ so that it equips $\A$ with a Poisson structure,
explicitly given by:
\begin{eqnarray}
\label{bracket_phi}
\PB_{\v}=\pp{\v}{z}\;\pp{}{x}\wedge\pp{}{y}+\pp{\v}{x}\;\pp{}{y}\wedge\pp{}{z}
+\pp{\v}{y}\;\pp{}{z}\wedge\pp{}{x}.
\end{eqnarray}

Our first purpose is to determine the Poisson
cohomology of this Poisson algebra $(\A,\PB_{\v})$, when $\v$ is a 
weight homogeneous polynomial with an isolated
singularity at the origin. 

We recall that the Poisson complex is constructed in the following way
(see \cite{zung} and \cite{plv} for details). First, the $k$-cochains of the Poisson 
complex of $(\A,\PB_{\v})$ are the skew-symmetric
$k$-derivations of $\A$ (i.e.\ the skew-symmetric $k$-linear maps $\A^k\rightarrow\A$
that are derivations in each of their arguments). We denote by
$\Vect^k(\A)$ the $\A$-module of all skew-symmetric
$k$-derivations of $\A$ and the elements of the $\A$-module
$
\Vect^*(\A)=\bigoplus_{k\in\N}\Vect^k(\A)
$ 
are called skew-symmetric multi-derivations of $\A$. By
convention, the $\A$-module of the $0$-derivations of $\A$ is
$\Vect^0(\A)=\A$.

The Poisson coboundary operator 
$\de^k_{\v}:\Vect^k(\A)\rightarrow\Vect^{k+1}(\A)$ is defined, 
for an element $Q\in\Vect^k(\A)$, by: 
\begin{eqn}{delta}
\lefteqn{\de^k_{\v}(Q)(f_0,\dots,f_k) :=
     \sum_{i=0}^k (-1)^i\pb{f_i,Q(f_0,\dots,\widehat{f_i},\dots,f_k)}_{\v}}\\
  &+&\displaystyle\sum_{0\leq i<j\leq k}(-1)^{i+j} 
Q(\pb{f_i,f_j}_{\v},f_0,\dots,\widehat{f_i}, 
        \dots,\widehat{f_j},\dots,f_k),
\end{eqn}%
where the symbol $\widehat{f_i}$ means that we omit the term $f_i$. It is
easy to see that $\de_{\v}^k(Q)$ is indeed a skew-symmetric
$(k+1)$-derivation while the fact that $\de_{\v}^{k+1}\circ\de_{\v}^k=0$ is an
easy consequence of the Jacobi identity (\ref{jacobi}).
The cohomology of this complex is called the Poisson cohomology of
$(\A,\PB_{\v})$.
We denote by $Z^k(\A,\v)$, respectively $B^k(\A,\v)$, the vector space
of all $k$-cocycles, respectively of all $k$-coboundaries, and we
denote by $H^k(\A,\v):=Z^k(\A,\v)/B^k(\A,\v)$, the $k$-th cohomology space.
As the space $H^0(\A,\v)$ is exactly the $\F$-vector space of the
Casimirs of $\PB_{\v}$ (i.e. the elements that belong to the
center of this bracket), we will also denote this space by $\Cas(\A,\v)$.
Notice that, if $\psi\in\Cas(\A,\v)$ , the operator $\de_{\v}$ commutes
with the multiplication by $\psi$. This implies that each of the
Poisson cohomology spaces $H^k(\A,\v)$ is a $\Cas(\A,\v)$-module.

In the case of the polynomial algebra $\A=\F[x,y,z]$, we have: 
\begin{eqnarray}
\label{isomultider}
\Vect^0(\A)\simeq\Vect^3(\A)\simeq \A, 
\qquad \Vect^1(\A)\simeq\Vect^2(\A)\simeq \A^3,
\end{eqnarray}
and $\Vect^k(\A)\simeq\ \{0\}$, for $k\geq 4$.
We choose these natural isomorphisms as follows:
$$
\renewcommand{\arraystretch}{1.3}
\begin{array}{ccc}
    \Vect^1(\A)&\longrightarrow& \A^3\\
       V &\longmapsto& (V[x],V[y],V[z]);
\end{array}\quad
\begin{array}{ccc}
    \Vect^2(\A)&\longrightarrow& \A^3\\
       V &\longmapsto& (V[y,z],V[z,x],V[x,y]);
\end{array}
$$
and
$\Vect^3(\A)\longrightarrow \A: V \longmapsto (V[x, y, z])$.

The elements of $\A^3$ are viewed as vector-valued functions on $\A$,
so we denote them with an arrow, like $\vec{f}\in\A^3$. Sometimes, it
will be important to distinguish $\A^3\simeq\Vect^1(\A)$ from 
$\A^3\simeq\Vect^2(\A)$; then we will rather write
$\vec{f}\in\Vect^1(\A)$ or $\vec{f}\in\Vect^2(\A)$. In $\A^3$,
let $\cdot$, $\times$ denote respectively the usual inner and cross
products, while $\vec{\nabla}$, $\vec{\nabla}\times$, $\div$ denote
respectively the gradient, the curl and the divergence operators.
For example, with these notations and the above isomorphisms, the
skew-symmetric biderivation $\PB_{\v}$ (defined in (\ref{bracket_phi})) 
is identified with the element
$\vn\v$ of $\A^3$.

Each of the Poisson coboundary operators $\de_{\v}^k$, given in (\ref{delta}), can
now be written in a compact form: 
\begin{eqn}{delta_i}
\de^0_{\v}(f) &=& \vn f\times\vn\v, \quad \hbox{ for }f\in\A\simeq\Vect^0(\A),\\
\de^1_{\v}(\vec{f}) &=& -\vn(\vec{f}\cdot\vn\v)+\div(\vec{f})\vn\v, 
\quad\hbox{ for } \vec{f}\in\A^3\simeq \Vect^1(\A),\\
\de^2_{\v}(\vec{f}) &=& -\vn\v\cdot(\vn\times\vec{f})=-\div(\vec{f}\times\vn\v),
\quad\hbox{ for } \vec{f}\in\A^3\simeq \Vect^2(\A),
\end{eqn}%
and the Poisson cohomology spaces of $(\A,\PB_{\v})$ take the
following forms
\begin{eqn*}[2.5]
 H^0(\A,\v)&=&\Cas(\A,\v)\simeq \{f\in\A \mid  \vn
                                  f\times\vn\v=\vec{0}\},\\
 H^1(\A,{\varphi})&\simeq& \frac{\displaystyle \{\vec{f}\in \A^3
                        \mid -\vn(\vec{f}\cdot\vn\v)
                           +\div(\vec{f})\vn\v=\vec{0}\}}
                     {\{\displaystyle
                        \vn f\times\vn\v\mid f\in\A\}},\\
H^2(\A,{\varphi})&\simeq& \frac{\displaystyle \{\vec{f}\in \A^3
                        \mid \vn\v\cdot(\vn\times\vec{f})=0\}}
                     {\{\displaystyle
                        -\vn(\vec{f}\cdot\vn\v)
                                +\div(\vec{f})\vn\v\mid
                        \vec{f}\in\A^3\}},\\
 H^3(\A,{\varphi})&\simeq&\frac{\displaystyle \A}{\displaystyle
  \{\vn\v\cdot(\vn\times \vec{f}) \mid \vec{f}\in\A^3\}}.
\end{eqn*}%

In order to compute these cohomology spaces, we will often use, 
for $\vec{f}, \vec{g}, \vec{h}\in\A^3$ and $f\in\A$, the following formulas,
well-known from vector calculus in~$\R^3$:
\begin{eqnarray}
\vn\times(f\vec{g}) &=& \vn f\times\vec{g}
      +f(\vn\times\vec{g}),\label{eq0}\\
\div(f \vec{g})&=&\vn f\cdot\vec{g}+f\div(\vec{g}),\label{divtimes3}\\
\div(\vec{f}\times\vec{g})&=& (\vn\times\vec{f})\cdot\vec{g}-
       \vec{f}\cdot(\vn\times\vec{g}).\label{divtimes2}
\end{eqnarray}
\subsection{Weight homogeneous multi-derivations}
\label{homder}
As we said, our results concern weight homogeneous Poisson
structures on $\A$. A  non-zero multi-derivation $P\in\Vect^*(\A)$ 
is said to be \textit{weight homogeneous} of (\textit{weighted}) degree 
$r\in\Z$, if there exist positive integers $\w_1, \w_2, \w_3\in\N^*$
(the \textit{weights} of the variables $x, y, z$), 
without a common divisor, such that $\L_{\vec{e}_{\w}}[P]=r P,$
where $\L_{\vec{e}_{\w}}$ is the Lie derivative with respect to the 
(weight homogeneous) Euler derivation 
$\vec{e}_{\w}=\w_1\,x\,\pp{}{x}+\w_2\,y\,\pp{}{y}+\w_3\,z\,\pp{}{z}$.
The degree of  a weight homogeneous multi-derivation $P\in\Vect^*(\A)$
is also denoted by $\w(P)\in\Z$. For $f\in\A$, it amounts to the usual
(weighted) degree of a polynomial. Notice that the degree of a
non-zero $k$-derivation may be negative for $k>0$. By convention, the
zero $k$-derivation is weight homogeneous of degree $-\infty$.

The Euler derivation $\vec{e}_{\w}$ is identified, with the
isomorphisms given in Section~\ref{ourstruct}, to the element 
$\vec{e}_{\w}=(\w_1\,x,\w_2\,y,\w_3\,z)\in\A^3$.
We denote by $|\w|$ the sum of the weights $\w_1+\w_2+\w_3$, so that 
$|\w|=\div(\vec{e}_{\w})$.
Euler's formula for a weight homogeneous $f\in\A$,
\begin{equation}\label{eulerhom1}
\vec{\nabla}f\cdot\vec{e}_{\w} = \w(f) f,
\end{equation}
then yields, using (\ref{divtimes3}):
\begin{equation}\label{eulerhom2}
\div(f\vec{e}_{\w}) = (\w(f)+|\w|) f.
\end{equation}

Fixing weights $\w_1, \w_2, \w_3\in\N^*$, it is clear that
$\A=\bigoplus_{i\in\N}\A_i$, where $\A_0=\F$ and for $i\in\N^*$,
$\A_i$ is the $\F$-vector space generated by all weight homogeneous polynomials
of degree $i$.
Denoting by $\Vect^k(\A)_i$ the $\F$-vector space
given by $\Vect^k(\A)_i:=\{P\in\Vect^k(\A)\mid\w(P)=i\}\cup\{0\}$, 
we have the following isomorphisms:
\begin{eqn}{corresdegree}
\Vect^0(\A)_i &\simeq& \A_i,\\
\Vect^1(\A)_i &\simeq& \A_{i+\w_1}\times\A_{i+\w_2}\times\A_{i+\w_3},\\
\Vect^2(\A)_i &\simeq& \A_{i+\w_2+\w_3}\times\A_{i+\w_1+\w_3}\times\A_{i+\w_1+\w_2},\\
\Vect^3(\A)_i &\simeq& \A_{i+\w_1+\w_2+\w_3}.
\end{eqn}%
Notice that even if $\Vect^1(\A)\simeq\Vect^2(\A)$ and
$\Vect^0(\A)\simeq\Vect^3(\A)$, these isomorphisms do not respect the
weight decompositions (\ref{corresdegree}).

One of our purposes is to determine the Poisson cohomology
of $(\A,\PB_{\v})$ when $\v\in\A$ is weight homogeneous with an
isolated singularity. The weight homogeneity of $\v$ 
will be essential for the computation of these spaces.
It implies indeed, among other things, that each of the coboundary 
operators $\de_{\v}^k$ is weight homogeneous of the same degree 
$N_{\w}:=\w(\v)-|\w|$, as can be seen from (\ref{delta_i}). That is to say, 
we have:
$$
P\in\Vect^k(\A)_i \Rightarrow \de_{\v}^k(P)\in\Vect^{k+1}(\A)_{i+N_{\w}}.
$$
If $P\in\Vect^k(\A)$ is a cocycle, then each of its 
weight homogeneous components will be a
cocycle. In the same way, if $P\in\Vect^k(\A)$ is a
coboundary then each of its weight homogeneous components will
be a coboundary. 
Moreover, if $P\in\Vect^k(\A)$ is a weight homogeneous coboundary, it
is the coboundary of a weight homogeneous element in $\Vect^{k-1}(\A)$.

\section{Isolated singularities and the Koszul complex}
\label{chapisokoszul}

In the next chapters, we will study the Poisson cohomology
associated to a weight homogeneous polynomial $\v\in\A=\F[x,y,z]$ (with $\Char(\F)=0$). 
As $\v$ will be supposed to have isolated singularities, we
will, in this part, recall some results about this notion, 
see \cite{sturm} and \cite{stan} for proofs.
 
Algebraically, we say that a weight homogeneous element $\v$ of $\F[x,y,z]$ 
has an \textit{isolated singularity} (at the origin) if  
\begin{eqnarray}
\label{A_sing}
\A_{sing}:=\F[x,y,z]/\ideal{\pp{\v}{x}, \pp{\v}{y}, \pp{\v}{z}}
\end{eqnarray}
is finite-dimensional, as a $\F$-vector space. The dimension of
$\A_{sing}$ is then called the Milnor number of the singular point. 
When $\F=\C$, this amounts, geometrically, to saying that the surface
$\Fe_{\v}:\{\v=0\}$ has a singular point only at the origin.
\begin{rem}
By definition, $\A_{sing}$ is exactly
the $\F$-algebra of regular functions of the affine variety 
$\Bigl\{\pp{\v}{x}= \pp{\v}{y}= \pp{\v}{z}=0\Bigr\}$ 
which is the singular locus of the
Poisson structure $\PB_{\v}$ (as can be seen from (\ref{bracket})). 
This algebra $\A_{sing}$ will play an important role in the
Poisson cohomology of the algebras $(\A,\PB_{\v})$ and $(\A_{\v},\PB_{\A_{\v}})$.
\end{rem}
Now, with the Cohen-Macaulay theorem, we will
see that, if $\v\in\A$ is a weight homogeneous
polynomial with an isolated singularity (what we will denote by
\textit{w.h.i.s.}), then the sequence of its partial derivatives $\pp{\v}{x},
\pp{\v}{y}, \pp{\v}{z}$ will be a regular sequence of $\A$. In order
to explain that, we first have to write down the definition of a 
homogeneous system of parameters of an algebra.
\begin{dfn}
Let $\A$ be an associative and commutative graded $\F$-algebra. A system of
homogeneous elements $F_1,\dots,F_d$ in $\A$, where $d$ is the Krull
dimension of $\A$, is called a \textit{homogeneous system of
parameters of $\A$} (\textit{h.s.o.p.})\ if $\A/\ideal{F_1,\dots,F_d}$
is a finite dimensional $\F$-vector space. 
\end{dfn}
For example, if we consider the $\F$-algebra $\A=\F[x,y,z]$, which is
graded by the weighted degree, we have a natural h.s.o.p.\ 
given by the system $x, y, z$.
Moreover, we have seen above that a weight homogeneous element
$\v\in\A$ has an isolated singularity (that is to say is w.h.i.s.)\ 
if and only if the three
partial derivatives $\pp{\v}{x}, \pp{\v}{y}, \pp{\v}{z}$ give a
h.s.o.p.\ of $\A$.

In order to understand the following theorem, that we will need, we still
have to give the definition of a regular sequence.
\begin{dfn}
A sequence $a_1, \dots, a_n$ in a commutative associative algebra $\A$
is said to be a \textit{$\A$-regular sequence} if
$\ideal{a_1,\dots,a_n}\not=\A$ and $a_i$ is  not a zero divisor
of $\A/\ideal{a_1,\dots,a_{i-1}}$ for $i=1,2,\dots,n$. 
\end{dfn}
For example, it is clear that the sequence $x, y, z$ is a regular
sequence in $\F[x,y,z]$. But, what about $\pp{\v}{x}, \pp{\v}{y},
\pp{\v}{z}$, when $\v$ is w.h.i.s.\ ?
\begin{thm}[Cohen-Macaulay]\label{cmthm}
Let $\A$ be a Noetherian graded $\F$-algebra. If $\A$ has a h.s.o.p.\ 
which is a regular sequence, then any h.s.o.p.\ in $\A$ is a
regular sequence. 
\end{thm}
Thus, when $\v\in\F[x,y,z]$ is w.h.i.s., then 
$\pp{\v}{x}, \pp{\v}{y}, \pp{\v}{z}$ is a regular sequence.
This is the key fact which leads to the following proposition, that
will play a fundamental role in our computations of Poisson cohomology, 
associated to a polynomial. 
%
\begin{prp}\label{diagram}
For any $\v\in\A$ the following diagram  
$$
\setlength{\dgARROWLENGTH}{0.6cm}
\begin{diagram}
\node{}\node{}\node{\F}\arrow{s,r}{}\node{\A}\arrow{s,r}{\vn}
              \node{\A^3}\arrow{s,r}{\vn\times}\\
\node{}\node{0}\arrow{e,t}{} 
         \node{\A} \arrow{s,r}{\vn} \arrow{e,t}{\vn\v}
             \node{\A^3} \arrow{e,t}{\times\vn\v} \arrow{s,r}{\vn\times}
             \node{\A^3} \arrow{e,t}{\cdot\vn\v} \arrow{s,r}{\div}
             \node{\A}\\
\node{}\node{\A} \arrow{e,t}{\vn\v} \arrow{s,r}{\vn}
            \node{\A^3} \arrow{e,t}{\times\vn\v} \arrow{s,r}{\vn\times}
            \node{\A^3} \arrow{e,t}{\cdot\vn\v} \arrow{s,r}{\div}
            \node{\A}\\
\node{\A}\arrow{e,t}{\vn\v}\node{\A^3}\arrow{e,t}{\times\vn\v}
           \node{\A^3}\arrow{e,t}{\cdot\vn\v}\node{\A}
\end{diagram}
$$
is commutative and has exact columns. If $\v$ is w.h.i.s.\ then the rows
of this diagram are also exact.
\end{prp}
\begin{rem}
\label{diaghom}
If $\v\in\A$ is weight homogeneous, then, as maps from
$\Vect^k(\A)$ to $\Vect^{k-1}(\A)$, each of the vertical arrows is
weight homogeneous of degree zero, while each of the horizontal arrows
is weight homogeneous of degree $\w(\v)$, the (weighted) degree of $\v$, leading to:
$$
\setlength{\dgARROWLENGTH}{0.6cm}
\begin{diagram}
\node{}\node{\Vect^3(\A)_{r}}\arrow{e,t}{\vn\v}\arrow{s,r}{\vn}
              \node{\Vect^2(\A)_{r+\w(\v)}}\arrow{s,r}{\vn\times}\\
\node{\Vect^3(\A)_{r-\w(\v)}} \arrow{s,r}{\vn} \arrow{e,t}{\vn\v}
             \node{\Vect^2(\A)_{r}} \arrow{e,t}{\times\vn\v} 
                \arrow{s,r}{\vn\times}
             \node{\Vect^1(\A)_{r+\w(\v)}} \arrow{e,t}{\cdot\vn\v} \arrow{s,r}{\div}
             \node{\Vect^0(\A)_{r+2\w(\v)}}\\
\node{\Vect^2(\A)_{r-\w(\v)}} \arrow{e,t}{\times\vn\v}
            \node{\Vect^1(\A)_{r}} \arrow{e,t}{\cdot\vn\v}
            \node{\Vect^0(\A)_{r+\w(\v)}}
\end{diagram}
$$
\end{rem}
\begin{proof}
Each column of this diagram is easily interpreted as the de Rham
complex of $\A$. The classical argument of exactness of the de Rham
complex of $C^{\infty}(\R^n)$ is easily adapted to the algebraic case: if
$\vec{f}=(f_1,f_2,f_3)\in\A^3$ is composed of three homogeneous 
polynomials of degree $d$ then $\div(\vec{f})=0$ implies 
that the first component of $\vn\times(\vec{f}\times\vec{e})$ is equal
to
$\left(\vn\times(\vec{f}\times\vec{e})\right)_1 
= 2 f_1+ \vn f_1\cdot\vec{e} -x\div(\vec{f})=(d+2)f_1$,
in view of Euler's Formula (\ref{eulerhom1}) 
($\vec{e}$ is the Euler derivation $(x,y,z)\in\A^3$, that is
to say $\vec{e}_{\w}$, with $\w_1=\w_2=\w_3=1$), so that 
$\vec{f}=\frac{1}{d+2}\vn\times(\vec{f}\times\vec{e})$.
Similarly, $\vn\times\vec{f}=\vec{0}$ implies that 
$\left(\vn(\vec{f}\cdot\vec{e})\right)_1 = f_1+\vn
f_1\cdot\vec{e}=(d+1)f_1$, that lieds to  
$\vec{f}=\frac{1}{d+1}\vn(\vec{f}\cdot\vec{e})$, according again to Euler's Formula.

Each of the rows of the diagram represents (part of) the so-called
Koszul complex. Let us prove that the Koszul complex, associated to 
$\v\in\A$ is exact, when $\v$ is w.h.i.s. If
$\vec{f}=(f_1, f_2, f_3)\in\A^3$ satisfies the equation 
$\vec{f}\times\vn\v=\vec{0}$, then we have three equalities 
like $f_1\pp{\v}{y}-f_2\pp{\v}{x}=0$. Since the partial derivatives of
$\v$ form a regular sequence, $\pp{\v}{y}$ is
not a zero divisor in $\A/\ideal{\pp{\v}{x}}$, so there exists
$\alpha\in\A$ such that $f_1=\alpha\pp{\v}{x}$ and then
$f_2=\alpha\pp{\v}{y}$. The other equations imply that
$f_3=\alpha\pp{\v}{z}$, that is to say $\vec{f}=\alpha\vn\v$. For the
second part of the exactitude of the Koszul complex, the reasoning is
exactly of the same kind.  
\end{proof}
\begin{rem}
\label{ifsquarefree}
If $\v\in\A$ is a weight
homogeneous polynomial without square factor then the first
part of the Koszul complex 
$
\A\stackrel{\vn\v}{\longrightarrow}\A^3\stackrel{\times\vn\v}{\longrightarrow}\A^3
$
is exact, but the second part 
$
\A^3\stackrel{\times\vn\v}{\longrightarrow}\A^3\stackrel{\cdot\vn\v}{\longrightarrow}\A
$
need not be exact if $\v$ is not w.h.i.s. For example, let
$\v=xyz\in\A$. The polynomial $\v$ is square free but the origin is not an
isolated singularity for $\v$. Then, the element $\vec{f}=(x,y,-2z)\in\A$
satisfies the equation $\vec{f}\cdot\vn\v=\vec{0}$ but, by an argument
of degree, there is no element $\vec{g}\in\A^3$ such that 
$\vec{f}=\vec{g}\times\vn\v$.
\end{rem}

We will often apply Proposition \ref{diagram} directly but sometimes,
we will use it in terms of the following corollary.
\begin{cor}\label{cor1}
Let $\v\in\A$ be w.h.i.s.\ and let $\vec{h}\in \A^3$. If
$(\vn\times\vec{h})\cdot \vn\v=0$ then there
exist $f,g\in \A$ such that $\vec{h}=\vn f+g\vn\v$.
\end{cor}
\begin{proof}
According to the diagram in Remark \ref{diaghom}, the operator 
$\vec{h}\mapsto (\vn\times\vec{h})\cdot \vn\v$, considered as a map
between $\Vect^2(\A)$ and $\Vect^0(\A)$, is a weight homogeneous
operator of degree $\w(\v)$. Therefore, 
it suffices to prove the result for an element 
$\vec{h}\in\Vect^2(\A)_r$, with $r\in\Z$.
If $(\vn\times\vec{h})\cdot \vn\v=0$ then, by Proposition
\ref{diagram}, there exists $\vec{k}\in\A^3$ such that 
$\vn\times\vec{h}=\vec{k}\times\vn\v$. In view of Remark
\ref{diaghom}, $\vec{k}$ can be chosen in $\Vect^2(\A)_{r-\w(\v)}$.
Summarizing, we have to prove that an equation of the type:
\begin{eqnarray}
\label{hyp_cor1}
\vn\times\vec{h}=\vec{k}\times\vn\v, 
\quad \vec{h}\in\Vect^2(\A)_r,\; \vec{k}\in\Vect^2(\A)_{r-\w(\v)}
\end{eqnarray}
implies that $\vec{h}=\vn f+g\vn\v$, with $f,g\in\A$.

We will do this by induction on $r\in\Z$, by proving the result directly for all 
$r<\w(\v)-\w^{[2]}$, with $\w^{[2]}:=\max\{\w_1+\w_2,\, \w_1+\w_3,\,
\w_2+\w_3\}$, where the integers $\w_1, \w_2, \w_3$ are the weights of the variables
$x,y,z$.

If $r<\w(\v)-\w^{[2]}$ then, according to the decompositions 
in (\ref{corresdegree}), $\Vect^2(\A)_{r-\w(\v)}=\{0\}$ so that the 
equality (\ref{hyp_cor1}) leads to
$\vn\times\vec{h}=\vec{0}$. Using Proposition \ref{diagram}, we 
obtain $\vec{h}=\vn f$, with $f\in\A$ as required.

Let $r'\geq \w(\v)-\w^{[2]}$ and assume that (\ref{hyp_cor1})
implies, for all $r<r'$, the existence of $f,g\in\A$ such that 
$\vec{h}=\vn f+g\vn\v$. 
Let us suppose that an element $\vec{l}\in\Vect^2(\A)_{r'}$ satisfies
an equation like in (\ref{hyp_cor1}), namely, suppose that
there exists $\vec{h}\in\Vect^2(\A)_{r'-\w(\v)}$ such that 
\begin{eqnarray}
\label{eq_cor1}
\vn\times\vec{l}=\vec{h}\times\vn\v.
\end{eqnarray}
Then, $\vec{h}$ satisfies (\ref{hyp_cor1}), with
$r=r'-\w(\v)$. Indeed, computing the divergence of both summands of
(\ref{eq_cor1}) gives 
$(\vn\times\vec{h})\cdot\vn\v=0$ and using
Proposition \ref{diagram} once again leads to the existence of 
$\vec{k}\in\Vect^2(\A)_{r'-2\w(\v)}$ such that we have
$\vn\times\vec{h}=\vec{k}\times\vn\v$.
By induction hypothesis, 
there exist $f,g\in\A$ such that $\vec{h}=\vn f+g\vn\v$.
Then, using Formula (\ref{eq0}), we obtain
$\vn\times\vec{l}=\vec{h}\times\vn\v=\vn f\times\vn\v= \vn\times(f\vn\v)$.

We can now conclude with Proposition \ref{diagram} that
there exists $f'\in \A$ such that $\vec{l}-f\vn\v=\vn f'$. Hence the result.
\end{proof}
\begin{rem}
\label{writing_Z2}
As $Z^2(\A,\v)=\{\vec{h}\in\A^3\mid (\vn\times\vec{h})\cdot
\vn\v=0\}$, Corollary \ref{cor1} leads to the equality
$$
Z^2(\A,\v)=\{\vn f+g\vn\v\mid f,g\in\A\}.
$$
This identity will be useful when we will determine $H^2(\A,\v)$ in
Section \ref{H2}.
\end{rem}
\section{Poisson cohomology associated to a weight homogeneous
  polynomial with an isolated singularity}
\label{computations}

Let us consider the polynomial algebra $\A=\F[x,y,z]$ ($\Char(\F)=0$), equipped with
the Poisson structure $\PB_{\v}$, where $\v\in\A$ is w.h.i.s.\ (weight
homogeneous polynomial with an isolated singularity).
We determine the Poisson cohomology spaces of the Poisson algebra
$(\A,\PB_{\v})$.

\begin{rem}\label{rem0}
If $\v\in\A$ is w.h.i.s.\ then $\w(\v)-\w_i>0$, for
$i=1,2,3$ (where $\w(\v)$ is still the (weighted) degree of $\v$ and 
$\w_1, \w_2, \w_3$ are
the weights of the variables $x,y,z$), and in particular, $\w(\v)>1$.
\end{rem}
\subsection{The space $H^0(\A,\v)$}
\label{H0}

A precise description of the $0$-th Poisson cohomology space,
which is also the algebra of the Casimirs, is given in the following proposition.
\begin{prp}\label{calculh0}
If $\v\in\A$ is w.h.i.s.\ then the zeroth Poisson cohomology space of 
$(\A,\PB_{\v})$ is given by 
$$
H^0(\A,\v)=\Cas(\A,\v)\simeq \bigoplus_{i\in\N} \F \v^i .
$$
\end{prp}
\begin{proof}
Let $f\in\A-\{0\}$ be a weight homogeneous $0$-cocycle, 
thus satisfying $\de_{\v}^0(f)=\vn f\times\vn\v=\vec{0}$.
Write $f$ as $f=h\v^r$, where $r\in\N$ and where $h\in\A-\{0\}$ is a
polynomial that is not divisible by $\v$. We have 
$\vn f= \v^r\vn h+rh\v^{r-1}\vn\v$, so $\vn h\times\vn\v=\vec{0}$.
Proposition \ref{diagram} implies the existence of $g\in\A$ such that 
$\vn h=g\vn\v$. Since $h$ and $\v$ are weight homogeneous and in view 
of Euler's Formula (\ref{eulerhom1}),
$$
\w(h)\, h =\vn h\cdot\vec{e}_{\w}=g\vn\v\cdot\vec{e}_{\w}=\w(\v)\, g\v,
$$
so $\w(h)=0$, as $h$ is not divisible by $\v$. Thus $h\in\F$ and 
$f=h\v^r\in\bigoplus_{i\in\N} \F \v^i$. 
Conversely, it is clear that 
$\de_{\v}^0(\v^r)=\vn(\v^r)\times\vn\v=\vec{0}$, for any $r\in\N$.
\end{proof}
\begin{rem}
According to Remark \ref{ifsquarefree}, if $\v\in\A$ is a weight
homogeneous polynomial without square factor but $\v$ is not necessarly
w.h.i.s., then the first
part of the Koszul complex 
is still exact, so Proposition \ref{calculh0} is also valid for this
more general class of polynomials.
However, if $\v$ has a square factor, the result is not true anymore. 
For example, if $\v=\psi^r$ with
$r\geq 2$ and $\psi\in\A$ a weight homogeneous polynomial without
square factor, then $H^0(\A,\v)\simeq H^0(\A,\psi)\simeq\bigoplus_{i\in\N}
\F \psi^i$ so that $H^0(\A,\v)\not\simeq\bigoplus_{i\in\N}
\F \v^i$.
\end{rem}
\subsection{The space $H^1(\A,\v)$}
\label{H1}

We first prove a result which will be useful to determine $H^1(\A,\v)$.
\begin{lma}
\label{ressyst}
Let $\v\in\A$ be w.h.i.s.\ and $\vec{g}\in\A^3$. Suppose that there exist
$r\in\N$ and $\alpha\in\F$ such that  
\begin{eqnarray}\label{system}
\left\lbrace
 \begin{array}{ccl}
  \vec{g}\cdot\vn\v&=&0,\\
  \div(\vec{g})&=&\alpha \v^r.
 \end{array}
\right.
\end{eqnarray}
Then $\alpha=0$ (equivalently $\div(\vec{g})=0$).
\end{lma}
\begin{proof}
According to Remark \ref{diaghom}, the operator 
$\vec{g}\mapsto(\vec{g}\cdot\vn\v,\div(\vec{g}))$ (from
$\A^3$ to $\A^2$) restricts for any $d\in\Z$ to an operator between
$\Vect^1(\A)_d$ and $\Vect^0(\A)_{d+\w(\v)}\times\Vect^0(\A)_d$.
Therefore it suffices to prove the lemma for an
element $\vec{g}\in\Vect^1(\A)_d$, with $d\in\Z$. Suppose that such an
element $\vec{g}$ satifies (\ref{system}), then, according to 
Proposition \ref{diagram}, the first 
equation implies that there exists $\vec{k}\in\Vect^2(\A)_{d-\w(\v)}$, 
such that $\vec{g}=\vec{k}\times\vn\v$. We will apply induction on $r\in\N$. 
First, if $r=0$, then, according to Formula (\ref{divtimes2}),  
$\alpha=\div(\vec{g})=\div(\vec{k}\times\vn\v)=(\vn\times\vec{k})\cdot\vn\v$,
so that $\alpha=0$, for degree reasons.

Assume now that for some fixed $r\geq 0$, any $\vec{g}$ that satisfies
(\ref{system}) is divergence free. Suppose that $\vec{h}\in\A^3$
satisfies $\vec{h}\cdot\vn\v=0$ and $\div(\vec{h})=\alpha'\v^{r+1}$,
for some $\alpha'\in\F$. Writing $\vec{h}=\vec{k}\times\vn\v$, the
Formulas (\ref{divtimes2}), (\ref{eulerhom1}) and
(\ref{eulerhom2}) show that $\vec{g}:=\vn\times\vec{k}-\frac{\alpha'}{\w(\v)}
\v^r \vec{e}_{\w}$ satisfies (\ref{system}), with 
$\alpha=-\alpha'(\w(\v)r+|\w|)/\w(\v)$, so that, by induction hypothesis, 
$0=\alpha=-\alpha'(\w(\v)r+|\w|)/\w(\v)$. It follows that $\alpha'=0$.
\end{proof}

Now, we can give the main result of this Section. We recall that
$|\w|$ is the sum of the weights of the three variables $x,y,z$.
\begin{prp}
\label{calculh1}
If $\v\in\A$ is w.h.i.s., then the first Poisson cohomology space of
$(\A,\PB_{\v})$ is a free module over $\Cas(\A,\v)$, given by:  
$$
H^1(\A,\v)\simeq\left\{
                \begin{array}{ccc}
                  \{0\} & \hbox{ if }&  \w(\v)\not= |\w|;\\
                  \Cas(\A,\v)\,\vec{e}_{\w}=\bigoplus\limits_{i\in\N}\F\v^i\,\vec{e}_{\w}
                  &\hbox{ if }& \w(\v)=|\w|.
                \end{array}
               \right.
$$
\end{prp}
\begin{proof}
Let $\vec{f}\in \Vect^1(\A)$ be a non zero element of $Z^1(\A,\v)$, that is to
say, $\vec{f}\in\A^3$ satisfies the equation: 
\begin{eqnarray}\label{eqcocycle}
\vn(\vec{f}\cdot\vn\v)=\div(\vec{f})\,\vn\v.
\end{eqnarray}
According to Remark \ref{diaghom}, we suppose that
$\vec{f}$ is weight homogeneous. Our purpose is to write 
$\vec{f}=\vn k \times\vn\v+ \frac{c}{\w(\v)} \v^{r} \vec{e}_{\w}
           \in B^1(\A,\v)+\bigoplus_{i\in\N}\F\v^i\,\vec{e}_{\w}$, 
where $c=0$ if $\w(\v)\not= |\w|$ and $c$ need not be $0$ otherwise.
Our proof will be divided in three parts.

$\quad 1.$ First, using cocycle condition (\ref{eqcocycle}), we find an
 element $\vec{g}\in\A^3$ which satisfies the equations (\ref{system}).
This equality implies indeed that 
$\de_{\v}^0(\vec{f}\cdot\vn\v)=\vn(\vec{f}\cdot\vn\v)\times\vn\v=\vec{0}$,
 so that the weight homogeneous element $\vec{f}\cdot\vn\v$ of $\A$ is a
Casimir. According to Proposition \ref{calculh0}, there exist
$c\in\F$ and $r\in\N$ such that $\vec{f}\cdot\vn\v = c \v^{r+1}$. Using
 Equation (\ref{eqcocycle}) once more, we obtain $\div(\vec{f})=c(r+1)\v^r$.
Letting $\vec{g} := \vec{f}- \frac{c}{\w(\v)}\v^r \vec{e}_{\w}$, 
Formulas (\ref{eulerhom1}) and
(\ref{eulerhom2}) imply that $\vec{g}$ satisfies (\ref{system}),
where $\alpha=c(1-\frac{|\w|}{\w(\v)})$.
Lemma \ref{ressyst} leads to 
$$
\left\lbrace
\begin{array}{l}
\div(\vec{g}) = 0,\; \vec{g}\cdot\vn\v =0,\\ 0=c\left(1-\frac{|\w|}{\w(\v)}\right).
\end{array}
\right.
$$

$\quad 2.$ Now, we will show that if $\vec{g}\in\A^3$ satisfies 
$\div(\vec{g}) = 0$ and $\vec{g}\cdot\vn\v =0$,
then $\vec{g}\in B^1(\A,\v)$.
Let $\vec{g}$ be a such element. As $\vec{g}\cdot\vn\v =0$,
Proposition~\ref{diagram} implies the existence of an element 
$\vec{h}\in\A^3$ such that 
$\vec{g}=\vec{h}\times\vn\v$. Moreover, we have
\begin{eqnarray*}
0=\div(\vec{g}) = \div(\vec{h}\times\vn\v)= (\vn\times\vec{h})\cdot\vn\v.
\end{eqnarray*}
Corollary \ref{cor1} leads now to the existence of elements
$k,l\in\A$ such that $\vec{h}=\vn k+l\vn\v$, so that 
$\vec{g}=\vn k \times\vn\v=\de_{\v}^0(k)\in B^1(\A,\v)$.

$\quad 3.$ The first two parts of this proof lead to the existence of
 $k\in\A$ and $c\in\F$ such that
\begin{eqnarray}\label{finB1}
 \left\lbrace
  \begin{array}{l}
   \vec{f} = \vn k \times\vn\v+ \frac{c}{\w(\v)} \v^r\vec{e}_{\w},\\
   0=c\left(1-\frac{|\w|}{\w(\v)}\right).
  \end{array}
 \right.
\end{eqnarray}
Now, we have to consider two cases:
$\w(\v)\not=|\w|$ and $\w(\v)=|\w|$.

$\bullet$ If $\w(\v)\not=|\w|$ then $c=0$ and 
$\vec{f}=\vn k\times\vn\v=\de_{\v}^0(k)\in B^1(\A,\v)$. Thus, when $\w(\v)\not=|\w|$, 
then $H^1(\A,\v)\simeq\{0\}$.

$\bullet$ Now, suppose that $\w(\v)=|\w|$, then (\ref{finB1}) leads to
$Z^1(\A,\v)\subseteq B^1(\A,\v)+\bigoplus_{i\in\N}\F\v^i\vec{e}_{\w}$. 
Conversely, for any $i\in\N$, Formulas (\ref{eulerhom1}) and (\ref{eulerhom2}) lead to  
$\de_{\v}^1(\v^i \vec{e}_{\w}) = (|\w|-\w(\v))\v^i \vn\v =0$. So that 
$$
Z^1(\A,\v)= B^1(\A,\v)+\bigoplus_{i\in\N}\F\v^i\vec{e}_{\w}.
$$
Let us show that this sum is a direct one.
It suffices to consider a weight homogeneous element 
$\alpha\v^i\vec{e}_{\w}\in B^1(\A,\v)$, $\alpha\in\F$, $i\in\N$.
It means that there exists $k\in\A$ such that 
$\alpha\v^i\vec{e}_{\w}=\vn k\times\vn\v$. Then (\ref{divtimes2}) 
and (\ref{eulerhom2}) lead to
$$
0=\div(\vn k\times\vn\v)=\div(\alpha\v^i\vec{e}_{\w})=\alpha|\w|(i+1)\v^i,
$$
therefore $\alpha=0$ and the sum 
$B^1(\A,\v)\oplus\bigoplus_{i\in\N}\F \v^i \vec{e}_{\w}$ is
direct. Thus, when $\w(\v)=|\w|$, then
$H^1(\A,\v)\simeq \bigoplus_{i\in\N}\F\v^i \vec{e}_{\w}$.
\end{proof}
\begin{rem}
We see that the case $\w(\v)=|\w|$ is particular. When $\v$ is homogeneous (i.e. weight
homogeneous with $\w_1=\w_2=\w_3=1$), it is the case where the degree
of $\v$ is three, that is
to say, where $\v$ is a cubic polynomial.
\end{rem}
\subsection{The space $H^3(\A,\v)$}
\label{H3}

Now, we give the third Poisson cohomology space of $(\A,\PB_{\v})$,
where $\v\in\A=\F[x,y,z]$ is w.h.i.s. Recall that, in this case,  
$$
\A_{sing}=\F[x,y,z]/\ideal{\pp{\v}{x}, \pp{\v}{y}, \pp{\v}{z}}
$$
is a finite dimensional $\F$-vector space, whose dimension is the 
Milnor number, denoted by $\mu$. Let $u_0=1,u_1,\dots,u_{\mu-1}$ be weight
homogeneous elements of $\A$, such that their images in $\A_{sing}$
give a $\F$-basis of
$\A_{sing}$. 

\begin{prp}
\label{calculh3}
If $\v\in\A=\F[x,y,z]$ is w.h.i.s.\ then the third cohomology space
$H^3(\A,\v)$ is the free $\Cas(\A,\v)$-module: 
\begin{eqnarray*}
H^3(\A,\v) \simeq \bigoplus_{j=0}^{\mu-1}\Cas(\A,\v)\; u_j 
\simeq \Cas(\A,\v)\otimes_{\F}\A_{sing}.
\end{eqnarray*}
\end{prp}
\begin{proof}
Let $f\in\A\simeq\Vect^3(\A)$ be a weight homogeneous polynomial of degree $d\in\N$.

$\quad 1.$ We first show that there exist $\vec{g}\in \A^3$, $N\in\N$ 
and elements $\lambda_{i,j}\in\F$, where $0\leq i\leq N$ and
$0\leq j \leq \mu-1$, such that: 
\begin{eqnarray}\label{elmntdeh3}
f=\vn\v\cdot(\vn\times\vec{g})+\sum_{i=0}^N 
              \sum_{j=0}^{\mu-1} \lambda_{i,j}\v^i u_j
\quad   \in B^3(\A,\v)+\sum_{\stackrel{k\in\N}{0\leq j\leq\mu-1}}\F\v^k u_j.
\end{eqnarray}
Let $\w^{[1]}:=\max(\w_1, \w_2, \w_3)$. We apply induction on $d$,
proving directly the result for $d\leq \w(\v)-\w^{[1]}$ (this is not
an empty case, as can be seen from Remark \ref{rem0}, for example,
it contains the case $f\in\F$).
By definition of the elements $u_0,\dots,u_{\mu-1}$, we have:
\begin{eqnarray}\label{defuj}
f=\vn\v\cdot\vec{l}+\sum_{j=0}^{\mu-1}\alpha_j u_j,
\end{eqnarray}
where $\vec{l}\in\Vect^1(\A)_{d-\w(\v)}$ and
$\alpha_0,\dots,\alpha_{\mu-1}\in\F$. 

If $d\leq \w(\v)-\w^{[1]}$ then the correspondences
(\ref{corresdegree}) imply that 
$\vec{l}$ is an element $(a,b,c)$ of $\F^3$ so that 
$f$ is indeed of the form (\ref{elmntdeh3}), with $\vec{g}=(bz, cx, ay)$, $N=0$
and $\lambda_{0,j}=\alpha_j$.

Now, suppose that $d> \w(\v)-\w^{[1]}$ and that any weight
homogeneous polynomial of degree at most $d-1$ is of the 
form (\ref{elmntdeh3}). Let us consider the decomposition
(\ref{defuj}) for $f$ of degree $d$. Proposition \ref{diagram} 
implies that there exists
$\vec{g}\in\A^3$ such that: 
\begin{eqnarray}
\label{forinduct}
\vec{l}-\frac{\div(\vec{l})}{d-\w(\v)+|\w|}\vec{e}_{\w}=\vn\times\vec{g},
\end{eqnarray}
since
$\div\biggl(\vec{l}-\frac{\div(\vec{l})}{d-\w(\v)+|\w|}\vec{e}_{\w}\biggr)=0$,
as follows from $\w(\div(\vec{l}))=d-\w(\v)$ and~(\ref{eulerhom2}).

Using the induction hypothesis on $\div(\vec{l})$, we conclude that
(\ref{defuj}), with $\vec{l}$ given by (\ref{forinduct}), is indeed of
the form (\ref{elmntdeh3}) (one uses that, according to Formula (\ref{eq0}),
$\v(\vn\times\vec{k})\cdot\vn\v=(\vn\times(\v\vec{k}))\cdot\vn\v$, for 
$\vec{k}\in\A^3$).
 
$\quad 2.$ So, we have already obtained that 
\begin{eqn}[2.1]{sumh3}
\A &=& \{\vn\v\cdot(\vn\times \vec{l}) \mid
  \vec{l}\in\A^3\}+ \displaystyle\sum_{j=0}^{\mu-1}\Cas(\A,\v)u_j\\
   &=& B^3(\A,\v)+ \displaystyle\sum_{j=0}^{\mu-1}\Cas(\A,\v)u_j.
\end{eqn}%
and it suffices to show that this sum is direct in $\A\simeq\Vect^3(\A)$.

We suppose the contrary.
This allows us to consider the smallest integer $N_0\in\N$
such that we have an equation of the form: 
\begin{eqnarray}\label{libertedsh3}
\sum_{i=N_0}^N \sum_{j=0}^{\mu-1} \lambda_{i,j}\v^i
u_j=\vn\v\cdot(\vn\times \vec{g})=-\de_{\v}^2(\vec{g}),
\end{eqnarray}
with $\vec{g}\in \A^3$, $N\geq N_0$ and $\lambda_{i,j}\in\F$ 
(for $N_0\leq i\leq N$ and $0\leq j \leq \mu-1$) and
$\lambda_{N_0,j_0}\not=0$, for some $0\leq j_0 \leq \mu-1$.
We will show that this hypothesis leads to a contradiction.

First, suppose that $N_0=0$, then the definition of the $u_j$, Euler's Formula
(\ref{eulerhom1}) and
(\ref{libertedsh3}) imply that $\lambda_{0,j}=0$ for all 
$0\leq j \leq \mu-1$, which contradicts the hypothesis $\lambda_{N_0,j_0}\not=0$.

So we suppose that $N_0>0$, using Euler's Formula
(\ref{eulerhom1}), the equation (\ref{libertedsh3}) can be written
as $\vn\v\cdot\Biggl(\sum_{i=N_0}^N \sum_{j=0}^{\mu-1} 
\frac{\lambda_{i,j}}{\w(\v)}\v^{i-1} u_j\vec{e}_{\w}\Biggr)
=\vn\v\cdot(\vn\times \vec{g})$.
Proposition \ref{diagram} implies that there exists
$\vec{h}\in\A^3$ such that: 
$$
\sum_{i=N_0}^N \sum_{j=0}^{\mu-1}
        \frac{\lambda_{i,j}}{\w(\v)}\v^{i-1} u_j\vec{e}_{\w}
=\vn\times \vec{g}+\vec{h}\times\vn\v.
$$
The divergence of both sides of this equality and Formula (\ref{eulerhom2}) give: 
$$
\sum_{i=N_1}^N \sum_{j=0}^{\mu-1} \lambda'_{i,j}\v^i u_j = 
(\vn\times\vec{h})\cdot\vn\v =-\de_{\v}^2(\vec{h}) ,
$$
where $\lambda'_{i,j}=\frac{\lambda_{i+1,j}}{\w(\v)}
(\w(\v)i+\w(u_j)+|\w|)$ and $N_1=N_0-1$.
So, we have obtained an equation of the form (\ref{libertedsh3}), with
$N_1<N_0$ and $\lambda'_{N_1,j_0}\not= 0$. This fact contradicts the
hypothesis and we conclude that the sum  (\ref{sumh3}) is direct.
The description of $H^3(\A,\v)$ follows.
\end{proof}
\subsection{The space $H^2(\A,\v)$}
\label{H2}

Finally, using Proposition \ref{calculh3} (and in fact the writing 
of $H^3(\A,\v)$), we obtain the second
Poisson cohomology space of the algebra $(\A,\PB_{\v})$, when
$\v\in\A=\F[x,y,z]$ is w.h.i.s. 

\begin{prp}
\label{calculh2}
If $\v\in\A=\F[x,y,z]$ is w.h.i.s.\ then the second Poisson
cohomology space of the algebra $(\A,\PB_{\v})$ is the $\Cas(\A,\v)$-module:
\begin{eqnarray*}
H^2(\A,\v)&\simeq& \bigoplus_{\stackrel{j=1}{\w(u_j)\not= \w(\v)-|\w|}}^{\mu-1} 
           \Cas(\A,\v)\vn u_j
                       \oplus
		   \bigoplus_{\stackrel{j=0}{\w(u_j)=\w(\v)-|\w|}}^{\mu-1} 
           \Cas(\A,\v)u_j\vn\v\\
       &&\qquad          \oplus 
                   \bigoplus_{\stackrel{j=1}{\w(u_j)=\w(\v)-|\w|}}^{\mu-1} 
           \F\vn u_j,
\end{eqnarray*}
where the first row gives the free part. 

In particular, we have: 
$H^2(\A,\v)\simeq\bigoplus_{j=1}^{\mu-1}
\Cas(\A,\v)\vn u_j$, if $\w(\v)<|\w|$ and 
$H^2(\A,\v)\simeq \bigoplus_{j=1}^{\mu-1} \Cas(\A,\v)\vn u_j \oplus
	\Cas(\A,\v)\vn\v$, when $\w(\v)=|\w|$.
\end{prp}
\begin{rem}
We see that the Poisson structure
$\PB_{\v}$ will be exact (that is to say a $2$-coboundary) if and only
if $\w(\v)\not=|\w|$. This fact comes from the equality 
$\de^1_{\v}(\vec{e}_{\w})=-(\w(\v)-|\w|)\vn\v$, a consequence of Formulas
(\ref{eulerhom1}) and (\ref{eulerhom2}).
\end{rem}
\begin{rem}
\label{magic_form}
Contrary to the other cohomology spaces, $H^2(\A,\v)$ is generally not
a free $\Cas(\A,\v)$-module. In fact, using Formulas (\ref{eulerhom1}) 
and (\ref{eulerhom2}), we get:
\begin{eqnarray}
\label{eq1}
\de_{\v}^1\left(\v^i u_j \vec{e}_{\w}\right) 
= \left(\w(u_j)-\w(\v)+|\w|\right)\v^i u_j\vn\v-\w(\v)\v^{i+1}\vn u_j.
\end{eqnarray}
This equality, which will be also useful later, explains that we 
have to distinguish, in the expression
of $H^2(\A,\v)$, the $u_j$ satisfying $\w(u_j)=\w(\v)-|\w|$ from the
other ones.
If $j$ is such that 
$\w(u_j)=\w(\v)-|\w|$ then (\ref{eq1}) yields that 
$\v^k\vn u_j\in B^2(\A,\v)$, for all $k\geq 1$, but this is not true
when $\w(u_j)\not=\w(\v)-|\w|$. This is the reason why
$H^2(\A,\v)$ is not always a free module over $\Cas(\A,\v)$.

Moreover, for all $j$ satisfying 
$\w(u_j)\not= \w(\v)-|\w|$, (\ref{eq1}) implies that $\v^i u_j \vn\v$,
$i\geq 0$, can be written as 
$c \v^{i+1}\vn u_j +
     \de_{\v}^1\left(c'\v^i u_j \vec{e}_{\w}\right)$,
with $c, c'\in\F-\{0\}$.
\end{rem}
\begin{proof}
First, let us show that: 
\begin{eqn}[2.1]{decomposition_H2}
Z^2(\A,\v)&\simeq& B^2(\A,\v)+\displaystyle\sum_{\stackrel{j=1}
     {\w(u_j)\not= \w(\v)-|\w|}}^{\mu-1} 
               \Cas(\A,\v)\vn u_j\\
          &+& \displaystyle\sum_{\stackrel{j=0}{\w(u_j)=\w(\v)-|\w|}}^{\mu-1} 
               \Cas(\A,\v)u_j\vn\v
           +  \displaystyle\sum_{\stackrel{j=1}{\w(u_j)=\w(\v)-|\w|}}^{\mu-1} \F\vn u_j.
\end{eqn}%

 Let $\vec{f}\in Z^2(\A,{\varphi})$. According to Remark \ref{writing_Z2}, 
there exists $g,h\in\A$ such that 
\begin{eqnarray}
\label{ecriturez2}
\vec{f}=\vn g+h\vn\v.
\end{eqnarray}
Moreover, Proposition \ref{calculh3} implies the existence of 
$\vec{g}_1,\vec{h}_1\in \A^3$, $N\in\N$ and of elements 
$\lambda_{i,j},\delta_{i,j}\in\F$, with 
$0\leq i\leq N$ and $0\leq j\leq\mu-1$, such that: 
\begin{eqnarray}\label{gandh}
g = \de_{\v}^2(\vec{g}_1)+\sum_{i=0}^N
\sum_{j=0}^{\mu-1} \lambda_{i,j}\v^i u_j,\quad
h = \de_{\v}^2(\vec{h}_1)+\sum_{i=0}^N
\sum_{j=0}^{\mu-1} \delta_{i,j}\v^i u_j,
\end{eqnarray}
while we have the $2$-coboundaries: 
\begin{eqnarray*}
\vn(\de_{\v}^2(\vec{g}_1)) &=& -\vn((\vn\times\vec{g}_1)\cdot\vn\v) 
           = \de_{\v}^1(\vn\times\vec{g}_1) \in B^2(\A,\v),\\
\de_{\v}^2(\vec{h}_1)\,\vn\v &=& \left((\vn\times\vec{h}_1)\cdot\vn\v\right)\vn\v 
           = \de_{\v}^1(\vec{h}_1\times\vn\v) \in B^2(\A,\v).
\end{eqnarray*}
Using this fact, (\ref{ecriturez2}) and (\ref{gandh}), we obtain 
\begin{eqnarray*}
\vec{f} \in B^2(\A,\v) + \sum_{j=1}^{\mu-1}  \Cas(\A,\v) \vn u_j
    + \sum_{j=0}^{\mu-1} \Cas(\A,\v) u_j \vn\v.
\end{eqnarray*}
Remark \ref{magic_form} then implies that $\vec{f}$ can be
decomposed as in the right hand side of~(\ref{decomposition_H2}). On
the other hand, all elements of the right hand side of 
(\ref{decomposition_H2}) are $2$-cocycles,
yielding equality in (\ref{decomposition_H2}). 
(Indeed, using Formula (\ref{eq0}), we have, for
all $f,g\in\A$, $\de_{\v}^2(\v\vn f)=-\vn\v\cdot(\vn\times(\v\vn f))=0$
and $\de_{\v}^2(g\vn\v)=-\vn\v\cdot(\vn\times(g\vn\v))=0$).
 
For the proof that the sum in (\ref{decomposition_H2}) is a direct
one, one uses the definition of the $u_j$ and applies
Propositions \ref{diagram}, \ref{calculh0} (expression of
$H^0(\A,\v)$) and \ref{calculh3} (writing of $H^3(\A,\v)$) as in the
proofs of Propositions \ref{calculh1} and \ref{calculh3}.
\end{proof}
\begin{rem}
Using Euler's
Formula (\ref{eulerhom1}) and the writings of the Poisson cohomology
spaces $H^1(\A,\v)$ and $H^2(\A,\v)$ given in Propositions
\ref{calculh1} and \ref{calculh2}, we can make the ring structure 
on the space 
$
H^*(\A,\v):=\bigoplus_{k=0}^3 H^k(\A,\v),
$
induced by the wedge product, explicit. One obtains, for example, that 
$\wedge : H^1(\A,\v)\times H^2(\A,\v) \longrightarrow H^3(\A,\v)$ is
surjective when $\w(\v)=|\w|$.
\end{rem}
\section{Poisson cohomology of the singular surface}
\label{surface}

In this chapter, we still consider an element $\v\in\A=\F[x,y,z]$ ($\Char(\F)=0$), which
is w.h.i.s.\ (weight homogeneous with an isolated singularity) and 
we restrict the Poisson structure
$\PB_{\v}$ to the singular surface $\Fe_{\v}:\{\v=0\}$ and compute
the cohomology of the Poisson algebra obtained.

\subsection{The Poisson complex of the singular surface $\Fe_{\v}$}
\label{complexsurf}
The algebra of regular functions on the surface $\Fe_{\v}$ is the
quotient algebra:
$$
\A_{\v}:= \frac{\F[x,y,z]}{\ideal{\v}}.
$$
Because $\v$ is a Casimir, $\ideal{\v}$ is a
Poisson ideal for $(\A,\PB_{\v})$ and the Poisson structure
$\PB_{\v}$ restricts naturally to $\Fe_{\v}$, that is to say goes down
to the quotient $\A_{\v}$. That leads to a Poisson structure on
$\A_{\v}$, denoted by $\PB_{\A_{\v}}$. Let us denote by $\pi$ the natural projection map 
$\A\rightarrow\A_{\v}$, then, for each $f,g\in\A$, we have 
$\pb{\pi(f),\pi(g)}_{\A_{\v}}=\pi\left(\pb{f,g}_{\v}\right)$ (that is
to say, $\pi$ is a Poisson morphism between $\A$ and $\A_{\v}$).
\begin{dfn}
\label{pi_related}
We say that $P\in\Vect^k(\A)$ and $Q\in\Vect^k(\A_{\v})$ are
\textit{$\pi$-related} and we write $Q=\pi_*(P)$ if
\begin{eqnarray}
\label{eq_pi_related}
\pi(P[f_1,\cdots,f_k])=Q[\pi(f_1),\cdots,\pi(f_k)],
\end{eqnarray}
for all $f_1,\cdots,f_k\in\A$.
\end{dfn}
In the following proposition, we give the Poisson cohomology spaces 
of the algebra $(\A_{\v},\PB_{\A_{\v}})$. That leads to consider the
skew-symmetric multi-derivations of the algebra
$\A_{\v}$ and the Poisson coboundary operators, associated to $\PB_{\A_{\v}}$. 
The previous definition will be useful in this discussion. 
By a slight abuse of notations we will, for an
element $\vec{f}=(f_1,f_2,f_3)\in\A^3$, denote by
$\pi(\vec{f})$, the element $(\pi(f_1),\pi(f_2),\pi(f_3))\in\A_{\v}^3$.
\begin{prp}
If $\v\in\A$ is w.h.i.s., the Poisson cohomology spaces of the algebra 
$(\A_{\v},\PB_{\A_{\v}})$, denoted by $H^k(\A_{\v})$, are given by:
\begin{eqnarray*}
\Cas(\A_{\v}) &=& H^0(\A_{\v})\simeq \left\{\pi(f)\in\A_{\v} \mid  
                \vn f\times\vn\v\in\ideal{\v}\right\},\\
             {}\\
H^1(\A_{\v}) &\simeq&  \frac{\displaystyle \left\{\pi\Bigl(\vec{f}\Bigr)\in \A_{\v}^3
                   \mid     \vec{f}\cdot\vn\v\in\ideal{\v} \hbox{ and }
            -\vn(\vec{f}\cdot\vn\v)  +\div(\vec{f})\,\vn\v\in\ideal{\v}\right\}}
        {\left\{\displaystyle  \pi\Bigl(\vn f\times\vn\v\Bigr)\mid  f\in\A\right\}},\\
             {}\\   
H^2(\A_{\v}) &\simeq& \frac{\displaystyle \left\{\pi\Bigl(\vec{f}\Bigr)\in\A_{\v}^3\mid 
                     \vec{f}\times\vn\v\in\ideal{\v}\right\}}
          {\left\{\displaystyle \pi \Bigl(-\vn(\vec{f}\cdot\vn\v)
           +\div(\vec{f})\,\vn\v\Bigr)\mid  \vec{f}\in\A^3;
                     \vec{f}\cdot\vn\v\in\ideal{\v}\right\}},
\end{eqnarray*}
and $H^3(\A_{\v})\simeq \{0\}$. 
\end{prp}
Subsequently, we denote by
$Z^k(\A_{\v})$ (respectively $B^k(\A_{\v})$) the space of all
$k$-cocycles (respectively $k$-coboundaries) of $\A_{\v}$.
\begin{proof}
We first have to determine the skew-symmetric multi-derivations of
$\A_{\v}$. Let us point out that any $P\in\Vect^k(\A)$ is $\pi$-related to a 
$Q\in\Vect^k(\A_{\v})$ if and only if $P[\v,f_2,\dots,f_k]\in\ideal{\v}$, for
all $f_2,\dots,f_k\in\A$. In this case, the equality
(\ref{eq_pi_related}) defines indeed an element $Q$ of $\Vect^k(\A_{\v})$, in
view of the skew-symmetry and the derivation properties of $P$.
Moreover, every $Q\in\Vect^k(\A_{\v})$ is obtained in this way.
Let us consider, for example, the case $k=1$. 

Let $Q\in\Vect^1(\A_{\v})$ and let us choose $\vec{f}=(f_1,f_2,f_3)\in\A^3$ such that 
$Q[\pi(x)]=\pi(f_1)$, $Q[\pi(y)]=\pi(f_2)$ and $Q[\pi(z)]=\pi(f_3)$.
Then, we get $Q=\pi_*(P)$, with
$P=f_1\pp{}{x}+f_2\pp{}{y}+f_3\pp{}{z}\in\Vect^1(\A)$ 
and $P[\v]=f_1\pp{\v}{x}+f_2\pp{\v}{y}+f_3\pp{\v}{z}=\vec{f}\cdot\vn\v\in\ideal{\v}$.

Conversely, each of 
$\pi(\vec{f})\in\A_{\v}^3$ satisfying the equation
$\vec{f}\cdot\vn\v\in\ideal{\v}$ 
gives an element of $\Vect^1(\A_{\v})$, defined by 
$\pi_*\left(f_1\pp{}{x}+f_2\pp{}{y}+f_3\pp{}{z}\right)$.
Thus, 
$$
\Vect^1(\A_{\v}) \simeq \{\pi(\vec{f})\in\A_{\v}^3 
              \mid \vec{f}\cdot\vn\v\in \ideal{\v}\}.
$$
With the same reasoning, we obtain 
$$
\Vect^2(\A_{\v}) \simeq \{\pi(\vec{f})\in\A_{\v}^3 \mid
               \vec{f}\times\vn\v\in\ideal{\v}\}.
$$
As it is clear that $\Vect^0(\A_{\v})\simeq\A_{\v}$ and 
$\Vect^k(\A_{\v})\simeq\{0\}$, for $k\geq
4$, let us now consider the space $\Vect^3(\A_{\v})$.
In the same way that above, we get 
$\Vect^3(\A_{\v})=\{\pi(f)\in\A_{\v}\mid
f\vn\v\in\ideal{\v}\}$. However, if $f\in\A$ satisfies
$f\vn\v=\v\,\vec{g}$, with $\vec{g}\in\A^3$, then we have 
$\vec{g}\times\vn\v=\vec{0}$ and
Proposition \ref{diagram} implies the existence of an element $h\in\A$
satifying $\vec{g}=h\vn\v$ so that $f=h\v\in\ideal{\v}$. That leads to
$
\Vect^3(\A_{\v})\simeq\{0\}.
$

Now, let us consider the Poisson coboundary operators of the Poisson
algebra $(\A_{\v},\PB_{\A_{\v}})$, denoted by $\de_{\A_{\v}}^k$. 
Using the definition of $\de_{\A_{\v}}^k$ (similarly as
(\ref{delta})), we obtain, for all $P\in\Vect^k(\A)$,
$\de_{\A_{\v}}^k(\pi_*(P))=\pi_*(\de_{\v}^k(P))$. 
That leads to:
\begin{eqn*}[1.5]
\de^0_{\A_{\v}}(\pi(f)) &=& \pi\left(\vn f\times\vn\v\right), 
           \quad \hbox{ for }\pi(f)\in\A_{\v}\simeq\Vect^0(\A_{\v}),\\
\de^1_{\A_{\v}}(\pi(\vec{f})) 
              &=& \pi\left(-\vn(\vec{f}\cdot\vn\v)+\div(\vec{f})\vn\v\right),\\
     &{}&\quad \quad\hbox{ for } \pi(\vec{f})\in\{\pi(\vec{g})\in\A_{\v}^3 
            \mid \vec{g}\cdot\vn\v\in
            \ideal{\v}\}\simeq\Vect^1(\A_{\v}),\\
\de^2_{\A_{\v}}(\pi(\vec{f})) &=& 0, 
  \quad\hbox{ for } \pi(\vec{f})\in\{\pi(\vec{g})\in\A_{\v}^3 \mid 
          \vec{g}\times\vn\v\in\ideal{\v}\}\simeq\Vect^2(\A_{\v}),
\end{eqn*}%
while the writing of the Poisson cohomology spaces follows.
\end{proof}
\subsection{The space $H^0(\A_{\v})$}
\label{h0surf}

In this Section, we consider still $\v\in\A$ w.h.i.s.\ and the Poisson
structure on $\A_{\v}$, denoted by $\PB_{\A_{\v}}$. We describe the
zeroth Poisson cohomology space, that is
to say the space of the Casimirs of $(\A_{\v},\PB_{\A_{\v}})$ in the
following Proposition.

\begin{prp}
\label{calculh0surf}
If $\v\in\A=\F[x,y,z]$ is w.h.i.s., the zeroth Poisson cohomology
space of the singular surface defined by this polynomial is given by 
$$
H^0(\A_{\v})=\Cas(\A_{\v})\simeq \F.
$$
\end{prp}
\begin{proof}
Let $f\in\A$ be a weight homogeneous polynomial such that $\pi(f)\in
H^0(\A_{\v})$. Then $\vn f\times\vn\v\in\ideal{\v}$
i.e., there exists $\vec{g}\in\A^3$ satifying 
$\vn f\times\vn\v=\v\, \vec{g}$. It follows that $\vec{g}\cdot\vn\v=0$ and
Proposition \ref{diagram} implies the existence of an element
$\vec{h}\in\A^3$ such that $\vec{g}=\vec{h}\times\vn\v$. Summing up, 
$(\vn f-\v\vec{h})\times\vn\v=\vec{0}$,
and we can apply Proposition \ref{diagram} again to obtain a $k\in\A$ 
satifying $\vn f=\v\vec{h}+k\,\vn\v$. Euler's
Formula (\ref{eulerhom1}) gives
$$
\w(f)\, f=\vn f\cdot\vec{e}_{\w}=\v(\vec{h}\cdot\vec{e}_{\w}+\w(\v)\,k).
$$
So, $f\in\ideal{\v}$ unless $\w(f)$,  the (weighted) degree of $f$, is
zero, thus $H^0(\A_{\v})\simeq\F$.
\end{proof}
\subsection{The space $H^1(\A_{\v})$}
\label{h1surf}

This section is devoted to the determination of the first Poisson
cohomology space of $(\A_{\v},\PB_{\A_{\v}})$, where
$\v\in\A=\F[x,y,z]$ is w.h.i.s. 
\begin{rem}
Using Proposition \ref{calculh0surf}, we can simplify the writing of
$Z^1(\A_{\v})$. Let indeed $\vec{f}\in\A^3$ be an element satisfying:
$-\vn(\vec{f}\cdot\vn\v)  +\div(\vec{f})\,\vn\v\in\ideal{\v}$.
Then $-\vn(\vec{f}\cdot\vn\v)\times\vn\v\in\ideal{\v}$, that is to say 
$\pi(\vec{f}\cdot\vn\v)\in H^0(\A_{\v})\simeq\F$, according to
Proposition \ref{calculh0surf}. For degree reasons, this leads to 
$\vec{f}\cdot\vn\v\in\ideal{\v}$. So, we can
simply write
$$
Z^1(\A_{\v})=\left\{\pi(\vec{f})\in \A_{\v}^3
            \mid -\vn(\vec{f}\cdot\vn\v)  +\div(\vec{f})\,\vn\v\in\ideal{\v}\right\}
$$
\end{rem}
Now, let us give the main result of this section (we recall that
$|\w|$ is the sum of the weights $\w_1,\w_2,\w_3$ of the variables
$x,y,z$ and that the family $\{u_j\}$ is an $\F$-basis of $\A_{sing}$ and is
defined in Section \ref{H3}).
\begin{prp}
\label{calculh1surf}
If $\v\in\A=\F[x,y,z]$ is w.h.i.s.\ then the first Poisson cohomology
space of the singular surface $\{\v=0\}$
is given by 
$$
H^1(\A_{\v})\simeq\bigoplus_{\stackrel{j=0}{\w(u_j)=\w(\v)-|\w|}}^{\mu-1}
            \F \pi(u_j\,\vec{e}_{\w}).
$$
In particular, if $\w(\v)<|\w|$ then $H^1(\A_{\v})\simeq\{0\}$.
\end{prp}
\begin{proof}
Let $\vec{f}\in\A^3$ satisfy $\pi\Bigl(\vec{f}\Bigr)\in
Z^1(\A_{\v})$, it means that there exists $\vec{k}\in\A^3$ 
satisfying $\de_{\v}^1(\vec{f})=\v\vec{k}$.
Then $0=\de_{\v}^2(\v\vec{k})=\v\,\de_{\v}^2(\vec{k})$,
because, as we said in Section~\ref{ourstruct}, the operator $\de_{\v}^2$ commutes
with the multiplication by $\v$. So $\v\,\vec{k}\in B^2(\A,\v)$ and
$\vec{k}\in Z^2(\A,\v)$. According to Proposition
\ref{calculh2},
\begin{eqnarray*}
\vec{k} &\in& B^2(\A,{\v}) \oplus
         \bigoplus_{\stackrel{j=1}{\w(u_j)\not=\w(\v)-|\w|}}^{\mu-1} 
                         \Cas(\A,{\v})\vn u_j\\
     & &    \oplus
    \bigoplus_{\stackrel{k=0}{\w(u_k)=\w(\v)-|\w|}}^{\mu-1} 
                         \Cas(\A,{\v})u_k\vn\v
          \oplus 
    \bigoplus_{\stackrel{l=1}{\w(u_l)=\w(\v)-|\w|}}^{\mu-1} 
                          \F\vn u_l\nonumber.
\end{eqnarray*}
Each of the first three summands is stable by multiplication by $\v$,
while Remark \ref{magic_form} gives
$$
\bigoplus_{\stackrel{l=1}{\w(u_l)=\w(\v)-|\w|}}^{\mu-1}
 \v\F\vn u_l\subset B^2(\A,{\v}).
$$
As a consequence, since $\v\vec{k}\in B^2(\A,\v)$, 
$$
\vec{k}\in B^2(\A,{\v})\oplus 
   \bigoplus_{\stackrel{l=1}{\w(u_l)=\w(\v)-|\w|}}^{\mu-1} \F\vn u_l.
$$
So there exist $\vec{h}\in\A^3$ and elements $\lambda_l\in\F$,
with $l$ satisfying $\w(u_l)=\w(\v)-|\w|$, such that
$$
\vec{k}=\de_{\v}^1(\vec{h})+\sum_{\stackrel{l=1}{\w(u_l)=\w(\v)-|\w|}}^{\mu-1}
                  \lambda_l \vn u_l.
$$
For all $1\leq l\leq \mu-1$ such that $\w(u_l)=\w(\v)-|\w|$, we have
$\v\vn u_l=-\de_{\v}^1\left(\frac{1}{\w(\v)} u_l\,\vec{e}_{\w}\right)$,
so that
$$
\de_{\v}^1(\vec{f})=\v\vec{k}=\de_{\v}^1\left(\v\vec{h}
-\sum_{\stackrel{l=1}{\w(u_l)=\w(\v)-|\w|}}^{\mu-1}\frac{\lambda_l}{\w(\v)}
u_l\, \vec{e}_{\w}\right).
$$
This implies
\begin{eqnarray}
\label{incomz1}
\vec{f}-\v\vec{h}+\sum_{\stackrel{l=1}{\w(u_l)=\w(\v)-|\w|}}^{\mu-1}
       \frac{\lambda_l}{\w(\v)} u_l\, \vec{e}_{\w}\;\in Z^1(\A,{\v}).
\end{eqnarray}

$\bullet$ If $\w(\v)\not=|\w|$, then Proposition
  \ref{calculh1} implies that (\ref{incomz1}) belongs to
  $B^1(\A,{\v})$, so that 
$$
\pi(\vec{f})\in \sum_{\stackrel{l=1}{\w(u_l)=\w(\v)-|\w|}}^{\mu-1}
       \F \pi(u_l\, \vec{e}_{\w}) + B^1(\A_{\v}).
$$

$\bullet$ If $\w(\v)=|\w|$ then
  (\ref{incomz1}) is simply the equation 
$\vec{f}-\v\vec{h}\in Z^1(\A,{\v})\simeq
  B^1(\A,{\v})+\Cas(\A,\v)\,\vec{e}_{\w}$, according to 
Proposition \ref{calculh1}. So, we have 
$\pi(\vec{f})\in \F\,\pi(\vec{e}_{\w})+B^1(\A_{\v})$.
As we have $\w(u_l)\geq 1$, if $1\leq l\leq \mu-1$, 
the result of both cases can be summarized as follows:
$$
Z^1(\A_{\v}) \subseteq B^1(\A_{\v})+\sum_{\stackrel{l=0}{\w(u_l)=\w(\v)-|\w|}}^{\mu-1}
                       \F \pi(u_l\,\vec{e}_{\w}).
$$
Euler's Formula (\ref{eulerhom1}) implies that
$\pi(u_l\,\vec{e}_{\w})\in Z^1(\A_{\v})$ ($\de_{\v}^1(u_l\,\vec{e}_{\w})\in \ideal{\v}$),
when $\w(u_l)=\w(\v)-|\w|$, so that the other inclusion holds too. 
It also allows us to show that the above sum
is a direct one. Hence the result about $H^1(\A_{\v})$.
\end{proof}
\subsection{The space $H^2(\A_{\v})$}
\label{h2surf}

We now compute the second Poisson cohomology space of $(\A_{\v},\PB_{\A_{\v}})$, where
$\v\in\A=\F[x,y,z]$ is w.h.i.s.
\begin{prp}
\label{calculh2surf}
If $\v\in\A=\F[x,y,z]$ is w.h.i.s.\ then $H^2(\A_{\v})$ is given by 
$$
H^2(\A_{\v})\simeq\bigoplus_{\stackrel{j=0}{\w(u_j)=\w(\v)-|\w|}}^{\mu-1}
\F \pi(u_j \vn\v).
$$
\end{prp}
\begin{rem}
It follows from Propositions \ref{calculh1surf} and \ref{calculh2surf}
that there is a natural isomorphism between $H^1(\A_{\v})$ and
$H^2(\A_{\v})$, that maps the element $u_j \,\vec{e}_{\w}$ (with 
$\w(u_j)=\w(\v)-|\w|$) to the element $u_j\,\vn\v$ of $H^2(\A_{\v})$.
\end{rem}
\begin{proof}
First, we show that the family $\left\{\pi(u_j\vn\v)\mid \w(u_j)=\w(\v)-|\w|\right\}$
generates the $\F$-vector space $H^2(\A_{\v})$.
Let $\vec{h}\in\A^3$ such that $\pi(\vec{h})\in
Z^2(\A_{\v})$, that is to say, such that there exists $\vec{g}\in\A^3$
satisfying $\vec{h}\times\vn\v=\v\,\vec{g}$.
According to Remark \ref{diaghom}, we may suppose
$\vec{h}\in\Vect^2(\A)_d$ and $\vec{g}\in\Vect^1(\A)_d$, with $d\in\Z$. 
Since $\vec{g}\cdot\vn\v=0$, Proposition \ref{diagram} 
implies that $\vec{g}=\vec{k}\times\vn\v$
and $\vec{h}=\v\vec{k}+f\vn\v$, with $f\in\Vect^3(\A)_{d-\w(\v)}$ and 
$\vec{k}\in\Vect^2(\A)_{d-\w(\v)}$.

If $d<\w(\v)-|\w|$ then $f=0$ and 
$\vec{h}\in\ideal{\v}$; otherwise $\pi(\vec{h})=\pi(f\vn\v)$, while, using Formulas
(\ref{eulerhom1}) and (\ref{eulerhom2}), we get
$\de_{\v}^1(f\vec{e}_{\w})=\left(d-2\w(\v)+2|\w|\right)f\vn\v-\w(\v)\,\v\vn
f$.
That leads, in the case $d\not=2(\w(\v)-|\w|)$, to 
$\pi(\vec{h})=\pi(f\vn\v)\in B^2(\A_{\v})$.
Therefore, let us suppose that $d=2(\w(\v)-|\w|)$, so that
$\w(f)=\w(\v)-|\w|$. For degree reasons, the projection map 
$\A\rightarrow\A_{sing}=\A/\ideal{\pp{\v}{x},\pp{\v}{y},\pp{\v}{z}}$
restricts to an injective map $\A_{\w(\v)-|\w|}\rightarrow\A_{sing}$, so that $f$ is
a $\F$-linear combination of the $u_j$ satisfying
$\w(u_j)=\w(\v)-|\w|$, that leads
to 
$$
\pi(\vec{h}) \in \sum_{\stackrel{j=0}{\w(u_j)=\w(\v)-|\w|}}^{\mu-1}\F \pi(u_j\vn\v),
$$
and for all $j$, $u_j\,\vn\v\in Z^2(\A_{\v})$.

It suffices now to show that this family is $\F$-free, modulo
$B^2(\A_{\v})$. It is empty if $\w(\v)<|\w|$, so we suppose
$\w(\v)\geq |\w|$. Let $\lambda_j$ be elements of $\F$ with $j$ such
that $\w(u_j)=\w(\v)-|\w|$ and let $\vec{l},\vec{\jmath}\in\A^3$ satisfying 
\begin{eqn}{cpresquefini}
\displaystyle\sum_{\stackrel{j=0}{\w(u_j) = \w(\v)-|\w|}}^{\mu-1}\lambda_j u_j\vn\v
           &=& -\vn(\vec{l}\cdot\vn\v)+\div(\vec{l})\vn\v+\v\vec{\jmath}\\
           &=& \de_{\v}^1(\vec{l})+\v\vec{\jmath},
\end{eqn}%
where the right hand side is an arbitrary
representative of an element of $B^2(\A_{\v})$. 
As the left hand side belongs to the space $\Vect^2(\A)_{2\w(\v)-2|\w|}$, 
we may suppose that $\vec{l}\in\Vect^1(\A)_{\w(\v)-|\w|}$
and $\vec{\jmath}\in\Vect^2(\A)_{\w(\v)-2|\w|}$.

The equation (\ref{cpresquefini}) implies 
$\vn(\vec{l}\cdot\vn\v)\times\vn\v\in\ideal{\v}$, so that 
$\pi\Bigl(\vec{l}\cdot\vn\v\Bigr)\in\Cas(\A_{\v})$. 
For degree reasons, Proposition \ref{calculh0surf} leads to the existence of $g\in\A$ 
of degree $\w(\v)-|\w|$ such that
$\vec{l}\cdot\vn\v=\v\, g=(g\vec{e}_{\w}\cdot\vn\v)/\w(\v)$. 
Then Proposition \ref{diagram} implies that $\w(\v)\,\vec{l}=g\vec{e}_{\w}$ 
and $\de_{\v}^1(\vec{l})=-\v\vn
g$, so that
\begin{eqnarray}\label{ca}
\sum_{\stackrel{j=0}{\w(u_j)=\w(\v)-|\w|}}^{\mu-1}\lambda_j u_j\vn\v=
           -\v\vn g+\v\vec{\jmath}=\v \vec{F},
\end{eqnarray}
where $\vec{F}=-\vn g+\vec{\jmath}\in\Vect^2(\A)_{\w(\v)-2|\w|}$. We get 
$\vec{F}\times\vn\v=\vec{0}$, but for degree reasons, Proposition
\ref{diagram} leads to $\vec{F}=\vec{0}$ so that, for all $j$, 
$\lambda_j=0$, since the family $\{u_j\}$ if
$\F$-free in $\A$.
\end{proof}

\section{Poisson homology associated to a weight homogeneous polynomial
  with an isolated singularity}
\label{homology}
In this last chapter, we consider the algebras 
$\A=\F[x,y,z]$ (with $\Char(\F)=0$) and $\A_{\v}=\A/\ideal{\v}$, where 
$\v\in\A$ is weight homogeneous with an 
isolated singularity (w.h.i.s.). These algebras are still
respectively equipped with the Poisson structures $\PB_{\v}$ and $\PB_{\A_{\v}}$. 
We use the Poisson cohomology of these Poisson algebras $(\A,\PB_{\v})$ and
$(\A_{\v},\PB_{\A_{\v}})$, given in the
previous chapters \ref{computations} and \ref{surface}, to determine
their Poisson homology.

\subsection{The Poisson homology of $\A$}
\label{homology_complex}
\subsubsection{Definitions}
\label{gen_hom}

We recall the construction of the Poisson homology complex 
associated to a Poisson algebra $(\B,\PB)$. First, the $k$-chains of this complex are
the so-called K\"ahler differential $k$-forms (see \cite{eis} for
details), whose space is denoted by $\Omega^k(\B)$. We recall that 
$\Omega^k(\B)=\wedge^k\Omega^1(\B)$ while 
$\Omega^*(\B):=\bigoplus_{k\in\N} \Omega^k(\B)$
is the $\B$-module of all (K\"ahler) differential forms, with, by
convention, $\Omega^0(\B)=\B$. We denote by $\diff$ the exterior
differential. The boundary operator, $\de_k:\Omega^k(\B)\rightarrow\Omega^{k-1}(\B)$,
called the Brylinsky or Koszul differential, is given by (see \cite{bryl}):
\begin{eqn}{Bdiff}
\lefteqn{\de_k(f_0\,\diff f_1\wedge\cdots\wedge\diff f_k)=
\sum_{i=1}^k(-1)^{i+1}\pb{f_0,f_i}\, 
    \diff f_1\wedge\cdots\wedge\widehat{\diff f_i}\wedge\dots\wedge \diff f_k}\\
  &+&\displaystyle\sum_{1\leq i<j\leq k} (-1)^{i+j} f_0\, \diff\pb{f_i,f_j}
        \wedge\diff f_1\wedge\cdots\wedge\widehat{\diff f_i}\wedge 
        \cdots\wedge\widehat{\diff f_j}\wedge\cdots\wedge \diff f_k,
\end{eqn}%
where the symbol $\widehat{\diff f_i}$ means that we omit the term
$\diff f_i$. It is easy to see that this operator satisfies
$\de_k\circ\de_{k+1}=0$. The homology of this complex 
is called the Poisson homology of
$(\B,\PB)$. 

The boundary operators of the algebras $(\A,\PB_{\v})$ and
$(\A_{\v},\PB_{\A_{\v}})$ are respectively
denoted by $\de_k^{\v}$ and $\de_k^{\A_{\v}}$, while the
Poisson homology spaces are denoted by $H_k(\A,\v)$ and
$H_k(\A_{\v})$. As for the Poisson cohomology, the 
boundary operator $\de_k$ commutes with the multiplication by a
Casimir, so that the Poisson homology spaces are modules over the
spaces of the Casimirs.

\subsubsection{The Poisson homology complex of $\A$}
\label{A_hom}

In the particular case of our polynomial algebra $\A=\F[x,y,z]$, it is clear that
$\Omega^*(\A)$ is the $\A$-module generated by the wedge products of
the $1$-differential forms $\diff x, \diff y, \diff z$ and that we
have $\Omega^i(\A)=\{0\}$, for all $i\geq 4$.
As for the multi-derivations of $\A$, we have the isomorphisms (with
the same choices as in Chapter \ref{ourstruct})
\begin{eqnarray}
\label{isodiffform}
\Omega^0(\A)\simeq\Omega^3(\A)\simeq \A, 
\qquad \Omega^1(\A)\simeq\Omega^2(\A)\simeq \A^3,
\end{eqnarray}
which allows us to use the same notations and formulas than in the previous
chapters, when we talk about differential forms. For example, the
$1$-differential form $\diff\v$ corresponds, with these notations, to
the element $\vn\v$ of $\A^3$ (as the biderivation $\PB_{\v}$).

\begin{prp}
\label{homology_A}
If $\v\in\A$ is w.h.i.s., the homology spaces of $(\A,\PB_{\v})$ are
given by:
$$
H_k(\A,\v)\simeq H^{3-k}(\A,\v), \hbox{ for all } k=0,1,2,3.
$$
\end{prp}
\begin{proof}
We have already seen in (\ref{isodiffform}) that
$\Omega^k(\A)\simeq\Vect^{3-k}(\A)$. In fact, for example, a
$1$-form $f\,\diff x\in\Omega^1(\A)$ corresponds to the biderivation 
$f\,\pp{}{y}\wedge\pp{}{z}\in\Vect^2(\A)$.
Moreover, under the previous identifications, we get easily 
$\de_k^{\v}=(-1)^k \de_{\v}^{3-k}$, that leads
to the result.
\end{proof}
\begin{rem}
There exists a more general result that gives, in certain cases, isomorphisms between
Poisson cohomology and homology spaces, using the modular class of a
Poisson algebra (see \cite{xu1} and \cite{hara} for details). 
\end{rem}
\subsection{The Poisson homology of $\A_{\v}$}
\label{surf_hom}
\subsubsection{The Poisson homology complex of $\A_{\v}$}

Now, let us determine the Poisson homology complex of the singular surface $\Fe_{\v}$.
For the quotient algebra $\A_{\v}=\F[x,y,z]/\ideal{\v}$, the space 
$\Omega^*(\A_{\v})$ is obtained by subjecting the $\A_{\v}$-module
generated by the wedge products of $\diff x, \diff y, \diff z$ to 
the relations $\v=0$, $\diff\v=0$ and $\diff\v\,\wedge\,\diff
x=0$,~etc.
We recall the natural surjective map $\pi:\A\rightarrow\A_{\v}$, which is a
Poisson morphism. This map induces another surjective map 
$\pi^{\sharp}:\Omega^k(\A)\rightarrow\Omega^k(\A_{\v})$ between the spaces
of all $k$-chains, which allows us to see the differential
$k$-forms of $\A_{\v}$ as images of differential $k$-forms of $\A$. Thus, as the
differential forms of $\A$ are identified to elements of $\A$ or $\A^3$, as
can be seen in (\ref{isodiffform}), we can write the spaces of all
differential $k$-forms of $\A_{\v}$ as quotients of $\A_{\v}$ and $\A_{\v}^3$
and then as quotients of $\A$ and $\A^3$.
We obtain, while $\Omega^0(\A_{\v}) \simeq \A_{\v}$,
\begin{eqnarray*}
\Omega^1(\A_{\v}) &\simeq& \frac{ \A_{\v}^3 }{ \{f\vn\v\mid f\in\A\} }
      \simeq \frac{\A^3}{ \{f\vn\v+\v\vec{g}\mid f\in\A,\,\vec{g}\in\A^3\} },\\
\Omega^2(\A_{\v}) &\simeq& \frac{ \A_{\v}^3 }
                          { \{\vn\v\times\vec{f}\mid\vec{f}\in\A^3\} }
      \simeq \frac{\A^3}{
	\{\vn\v\times\vec{f}+\v\vec{g}\mid\vec{f},\vec{g}\in\A^3\} },\\
\Omega^3(\A_{\v}) &\simeq& \frac{ \A_{\v} }
    { \{\vn\v\cdot\vec{f}\mid\vec{f}\in\A^3\} }
    \simeq \frac{\A}{\displaystyle \ideal{\pp{\v}{x},\pp{\v}{y},\pp{\v}{z}}}=\A_{sing}.
\end{eqnarray*}
\begin{rem}
Unlike for $\A$, there is no isomorphisms between
the spaces of skew-symmetric multi-derivations and differential
forms on $\A_{\v}$. For example, 
$\Omega^0(\A_{\v})\simeq\A_{\v}$ while 
$\Vect^3(\A_{\v})\simeq\{0\}$ and $\Vect^2(\A_{\v})\subseteq \A_{\v}^3$.
Observe also that $\Omega^3(\A_{\v})\not\simeq \{0\}$, although $\Fe_{\v}$ is an
affine variety of dimension two. 
\end{rem}
In view of Definition (\ref{Bdiff}), the operator $\de_k^{\v}$
induces an operator
$\Omega^k(\A_{\v})\rightarrow\Omega^{k-1}(\A_{\v})$, that is exactly $\de_k^{\A_{\v}}$,
so that the Poisson homology spaces of $\A_{\v}$ are given by
\begin{eqnarray*}
H_0(\A_{\v}) &\simeq& \frac{\displaystyle \A}
         {\displaystyle
         \{\vn\v\cdot(\vn\times\vec{f})+\v\,g\mid g\in\A,\vec{f}\in\A^3\}},\\
&&\\
H_1(\A_{\v}) &\simeq& 
\frac{\displaystyle \{\vec{f}\in\A^3\mid\vn\v\cdot(\vn\times\vec{f})\in\ideal{\v}\} }
{\displaystyle
         \{-\vn(\vec{f}\cdot\vn\v)+\div(\vec{f})\,\vn\v+g\vn\v+\v\vec{h}\mid
         g\in\A, \vec{f},\vec{h}\in\A^3\}},\\
&&\\
H_2(\A_{\v}) &\simeq& 
\frac{\displaystyle \{\vec{f}\in\A^3\mid
   -\vn(\vec{f}\cdot\vn\v)+\div(\vec{f})\,\vn\v
   \in{\mathcal I}_{\v}\}}
{\displaystyle \{\vn\v\times\vec{h}+\v\vec{k}\mid\vec{h},
         \vec{k}\in\A^3\}},\\
&&\qquad\qquad\qquad\hbox{ where } {\mathcal I}_{\v}:=\{f\vn\v+\v\vec{g}\mid f\in\A,
         \vec{g}\in\A^3 \} ,\\
H_3(\A_{\v}) &\simeq& \A_{sing}.
\end{eqnarray*}
\begin{rem}
In view of the writing of the Poisson homology groups of $\A$ and
$\A_{\v}$, we can describe explicitly the map induced by $\pi$ between
these groups. In fact, this map is exactly the reduction modulo $\v$
between the spaces $H_k(\A)$ and $H_k(\A_{\v})$, for $k\not=1$, and it is
the reduction modulo ${\mathcal I}_{\v}$, for $k=1$. This phenomenon will
be illustrated in the determination of the Poisson homology groups of $\A_{\v}$.
\end{rem}

\subsubsection{The Poisson homology spaces of the singular surface $\Fe_{\v}$}
\label{surf_hom_spaces}

In this Section, $\v\in\A=\F[x,y,z]$ is still w.h.i.s.\ and we determine these spaces.

\begin{prp}
\label{homology_surf}
If $\v\in\A$ is w.h.i.s.\ then the homology spaces of the singular
surface are given by:
\begin{eqnarray*}
H_0(\A_{\v}) &\simeq& \bigoplus_{j=0}^{\mu-1}\,\F u_j \simeq \A_{sing},\qquad
H_1(\A_{\v}) \simeq \bigoplus_{j=1}^{\mu-1}\,\F\vn u_j,\\
&&\quad H_2(\A_{\v}) \simeq \bigoplus_{j=0}^{\mu-1}\,\F u_j\,\vec{e}_{\w}\simeq\A_{sing}.
\end{eqnarray*}
\end{prp}
\begin{rem}
The fact that $H_0(\A_{\v})\simeq \A_{sing}$ was already proved by J. Alev
and T. Lambre, with other methods, in \cite{A_L}. Their result is more
general as they only suppose that $\v$ is a weight homogeneous polynomial,
not necessarly with an isolated singularity. 
\end{rem}
\begin{rem}
The multiplication by $\vec{e}_{\w}$ gives a natural isomorphism between
$H_0(\A_{\v})$ and $H_2(\A_{\v})$, while the operator of gradient $\vn$
gives a surjective map from $H_0(\A_{\v})$ to $H_1(\A_{\v})$.
\end{rem}
\begin{proof}
$\quad 1.$ We first determine $H_0(\A_{\v})$. 
According to Proposition \ref{calculh3} (i.e., the writing of
$H^3(\A,\v)$), we have: 
\begin{eqnarray*}
\A &=& \{\vn\v\cdot(\vn\times\vec{f})\mid\vec{f}\in\A^3\}
       +\sum_{\stackrel{j=0}{i\in\N}}^{\mu-1}\F\v^i u_j,\\
   &=& \{\vn\v\cdot(\vn\times\vec{f})+\v g\mid g\in\A, \vec{f}\in\A^3\}
       +\bigoplus_{j=0}^{\mu-1}\F\, u_j.
\end{eqnarray*}
Moreover this last sum is a direct one, as follows from the definition
of the $u_j$ (in Section \ref{H3}) and the inclusion
$\{\vn\v\cdot(\vn\times\vec{f})+\v g\mid g\in\A,
\vec{f}\in\A^3\}\subseteq \ideal{\pp{\v}{x}, \pp{\v}{y},
  \pp{\v}{z}}$, easily obtained with Euler's Formula
(\ref{eulerhom1}). That leads to
$H_0(\A_{\v})\simeq\bigoplus_{j=0}^{\mu-1}\F\,u_j$.\\

$\quad 2.$ Now, we use the result we obtained for $H^2(\A,\v)$ to determine the
first Poisson homology space of $\A_{\v}$.
Let $\vec{f}\in\A^3$ satisfying
$\vn\v\cdot(\vn\times\vec{f})\in\ideal{\v}$, thus,
there exists $g\in\A$ with $-\de_{\v}^2(\vec{f})=\vn\v\cdot(\vn\times\vec{f})=\v\,g$.

According to Proposition \ref{calculh3}, 
$g\in B^3(\A,\v) \oplus\bigoplus_{j=0}^{\mu-1}\Cas(\A,\v)\,u_j$.
As both of the summands of this sum are stable by multiplication by
$\v$ and because $\v g\in B^3(\A,\v)$, we have $g\in B^3(\A,\v)$, i.e. there
exists $\vec{k}\in\A^3$ satisfying $g=\vn\v\cdot(\vn\times\vec{k})$. 
Thus, $\vec{f}-\v\vec{k}\in Z^2(\A,\v)$ together with Proposition
\ref{calculh2} imply that 
$$
\vec{f}\in \sum_{j=1}^{\mu-1} \F\,\vn u_j+\{\de_{\v}^1(\vec{l})+g\vn\v+\v\vec{h}\mid
         g\in\A, \vec{l},\vec{h}\in\A^3\},
$$
so that $\{\vn u_j\mid 1\leq j\leq\mu-1\}$ generates the
$\F$-vector space $H_1(\A_{\v})$ and it suffices to prove that $\vn
u_1,\dots,\vn u_{\mu-1}$ are linearly independent elements of $H_1(\A_{\v})$.
Assume therefore that there exist elements $\lambda_j$ of $\F$ ($1\leq
j\leq\mu-1$), $\vec{k},\vec{l}\in\A^3$ and $g\in\A$ such that
$$
\sum_{j=1}^{\mu-1} \lambda_j \vn u_j=
-\vn(\vec{l}\cdot\vn\v)+\div(\vec{l})\,\vn\v+g\vn\v+\v\vec{h}.
$$
Then, as the $u_j$ are weight homogeneous, Euler's Formula (\ref{eulerhom1}) leads to 
$$
\sum_{j=1}^{\mu-1} \lambda_j \w(u_j)\,u_j\in\ideal{\pp{\v}{x},
  \pp{\v}{y}, \pp{\v}{z}}
$$ 
and the definition of the $u_j$ implies
$\lambda_j=0$, for $1\leq j\leq\mu-1$.\\

$\quad 3.$ Finally, we compute the second Poisson homology space of $\A_{\v}$.
Let $\vec{f}\in\A^3$ satisfying $\de^1(\vec{f})\in{\mathcal I}_{\v}$,
i.e. there exist $l\in\A$, $\vec{g}\in\A^3$ such that 
$\de^1(\vec{f})=l\vn\v+\v\,\vec{g}$.

$\bullet$ Let us study the term $\v\,\vec{g}$. We first point out 
that $l\vn\v\in Z^2(\A,\v)$, so that 
$\v\,\vec{g}=\de^1(\vec{f})-l\vn\v\in Z^2(\A,\v)$.
Using Proposition \ref{calculh2}, Formula (\ref{eq1})
and the fact that $\de^1_{\v}$ commutes with $\v$, we obtain the
existence of $\vec{h}\in\A^3$ and $c_j\in\F$, such that:
\begin{eqn}[2.5]{Z2*phi}
\v\,\vec{g} &\in& \de_{\v}^1\left(\v\,\vec{h}
     +\displaystyle\sum_{\stackrel{j=1}{\w(u_j)=\w(\v)-|\w|}}
         ^{\mu-1}c_j u_j\,\vec{e}_{\w}\right)\\
     &+& \displaystyle\bigoplus_{\stackrel{j=1}{\w(u_j)\not= \w(\v)-|\w|}}^{\mu-1} 
           \Cas(\A,\v)\vn u_j     \oplus
      \displaystyle\bigoplus_{\stackrel{j=0}{\w(u_j)=\w(\v)-|\w|}}^{\mu-1} 
           \Cas(\A,\v)u_j\vn\v,
\end{eqn}%

$\bullet$ Next, we consider the term $l\vn\v$. According to Proposition
  \ref{calculh3}, there exists $\vec{k}\in\A^3$ such that 
$l\in \de_{\v}^2(\vec{k})+\Cas(\A,\v)\otimes_{\F}\A_{sing}$. The
  equality
  $\de_{\v}^2(\vec{k})\,\vn\v=\de_{\v}^1(\vec{k}\times\vn\v)$ and
  Formula (\ref{eq1}) lead to:
\begin{eqn}[2.5]{A*grad_phi}
l\vn\v &\in& \de_{\v}^1\left(\vec{k}\times\vn\v 
         + \displaystyle\sum_{\stackrel{j=0}{\w(u_j)\not=\w(\v)-|\w|}}^{\mu-1} 
            \mathcal{C}_j\,  u_j\,\vec{e}_{\w}\right)\\
         &+& \displaystyle\bigoplus_{\stackrel{j=1}{\w(u_j)\not= \w(\v)-|\w|}}^{\mu-1} 
           \Cas(\A,\v)\vn u_j     \oplus
     \displaystyle \bigoplus_{\stackrel{j=0}{\w(u_j)=\w(\v)-|\w|}}^{\mu-1} 
           \Cas(\A,\v)u_j\vn\v,
\end{eqn}%
where $\mathcal{C}_j\in\Cas(\A,\v)$. \\

The equalities (\ref{Z2*phi}) and (\ref{A*grad_phi}) give:
\begin{eqnarray*}
\lefteqn{\de_{\v}^1\left(\vec{f}-\v\,\vec{h}
 -\sum_{\stackrel{j=1}{\w(u_j)=\w(\v)-|\w|}}^{\mu-1}c_j u_j\,\vec{e}_{\w}
 -\vec{k}\times\vn\v 
         - \sum_{\stackrel{j=0}{\w(u_j)\not=\w(\v)-|\w|}}^{\mu-1} 
       \mathcal{C}_j\,u_j\,\vec{e}_{\w}\right)}\\
&\in& \bigoplus_{\stackrel{j=1}{\w(u_j)\not= \w(\v)-|\w|}}^{\mu-1} 
           \Cas(\A,\v)\vn u_j     \oplus
      \bigoplus_{\stackrel{j=0}{\w(u_j)=\w(\v)-|\w|}}^{\mu-1} 
           \Cas(\A,\v)u_j\vn\v.
\end{eqnarray*}
Using Proposition \ref{calculh2} once more, we obtain
$$
\vec{f}-\v\,\vec{h}
 -\sum_{\stackrel{j=1}{\w(u_j)=\w(\v)-|\w|}}^{\mu-1}c_j u_j\,\vec{e}_{\w}
 -\vec{k}\times\vn\v 
         - \sum_{\stackrel{j=0}{\w(u_j)\not=\w(\v)-|\w|}}^{\mu-1} 
       \mathcal{C}_j\,u_j\,\vec{e}_{\w} \in Z^1(\A,\v).
$$
It suffices now to use Proposition \ref{calculh1} to conclude that
$$
\vec{f}\in \sum_{j=0}^{\mu-1} \F\,u_j\,\vec{e}_{\w} + 
\{\vn\v\times\vec{k}+\v\vec{h}\mid\vec{h},\vec{k}\in\A^3\}.
$$
Finally, using Euler's Formula (\ref{eulerhom1}) and the definition of the
$u_j$, it is easy the see that this sum is a direct one in
$\A^3$. Hence the result for $H_2(\A_{\v})$.
\end{proof}
\bibliographystyle{plain}
\bibliography{ref}
%


\end{document}